\documentclass[12pt]{amsart}
\usepackage[usenames,dvipsnames]{xcolor}
\usepackage{amsmath,amssymb,amsthm,graphicx,mathrsfs,url,latexsym,enumerate}
\usepackage{a4wide}
\usepackage{ifthen}
\usepackage{verbatim}
\usepackage{fancyhdr}
\usepackage{url}
\usepackage{bbm}
\usepackage[sans]{dsfont}
\usepackage{transparent}
\usepackage{amsbsy}
\usepackage[colorlinks=true,linkcolor=Red,citecolor=Green]{hyperref}
\usepackage{wrapfig}
\usepackage{tikz}

\definecolor{bleu}{RGB}{27,88,145}
\definecolor{mauve}{RGB}{138,20,79}

\renewcommand{\Re}{\operatorname{Re}}

\newcommand{\Id}{\operatorname{Id}}

\newcommand{\supp}{\operatorname{supp}}

\renewcommand{\span}{\operatorname{span}}
\newcommand{\C}{\mathbb C}

\newcommand{\N}{\mathbb N}
\newcommand{\R}{\mathbb R}

\newcommand{\Hess}{\operatorname{Hess}}

\newcommand{\Ran}{\operatorname{Ran}}
\newcommand{\Ker}{\operatorname{Ker}}

\renewcommand{\div}{\operatorname{div}}
\def\<{\langle}
\def\>{\rangle}
\newcommand{\bp}{{\it Proof. }}
\newcommand{\ep}{\hfill $\square$}
\newcommand{\be}{\begin{equation}}
\newcommand{\ee}{\end{equation}}
\newcommand{\bes}{\begin{equation*}}
\newcommand{\ees}{\end{equation*}}

\numberwithin{equation}{section}
\numberwithin{figure}{section}

\def\m{\mathbf{m}}
\def\s{\mathbf{s}}

\everymath{\displaystyle}
\renewcommand{\Pr}{\Pi_\rho}
\DeclareMathOperator{\p}{\partial}
\DeclareMathOperator{\eps}{\varepsilon}
\DeclareMathOperator{\phii}{\varphi}

\newcommand{\dz}{\ dz}

\newcommand{\dv}{\ dv}
\newcommand{\vvert}[1]{\left\lVert#1\right\rVert}

\newtheorem{theorem}{Theorem}
\newtheorem{defin}{Definition}[section]

\newtheorem{lemma}[defin]{Lemma}

\newtheorem{proposition}[defin]{Proposition}
\newtheorem{remark}[defin]{Remark}

\newtheorem{assumption}{Assumption}

\def\bbb{{\mathcal B}}\def\ccc{{\mathcal C}}\def\ddd{{\mathcal D}}
 
 \def\lll{{\mathcal L}}
 \def\ooo{{\mathcal O}}\def\ppp{{\mathcal P}}

\def\uuu{{\mathcal U}}\def\www{{\mathcal W}}

\begin{document}
\title{Sharp spectral gap of adaptive Langevin dynamics}
\maketitle

\begin{center}\textsc{Loïs Delande}\end{center}

\begin{abstract}
    We consider a degenerated Fokker-Planck type differential operator associated to an adaptive Langevin dynamic. We prove Eyring–Kramers formulas for the bottom of the spectrum of this operator in the low temperature regime. The main ingredients are resolvent estimates obtained via hypocoercive techniques and the construction of sharp Gaussian quasimodes through an adaptation of the WKB method.
\end{abstract}

\section{Introduction}

\subsection{Motivations.}
In order to describe the dynamic that rules the evolution of a molecular system at temperature of order $h>0$, the following homogeneous Langevin process is widely used
\be\label{eq:SDE0}
dX_t = \xi(X_t)dt + \sqrt{2h}\sigma(X_t)dB_t,
\ee
where $(X_t)_{t\geq0}$ gives the positions of the particles, the vector field $\xi$ is the drift coefficient, the matrix field $\sigma$ is the diffusion coefficient and $(B_t)_{t\geq0}$ denotes a $d$-dimensional Brownian motion. In the low temperature regime, i.e. when $h\to0$, we observe a metastable behaviour of the solution of \eqref{eq:SDE0}, this can be obtained via the study of the exit problem for this SDE. Considering an open set $\Omega$ and a point $x\in\Omega$ the question is to know where and when does the process exit $\Omega$ having set $X_0=x$. This problem has been intensively studied in \cite{FrWe84} or \cite{Da83} for some pioneer work, we also refer to \cite{DiGLeLLePNec18}, \cite{Ne21} and \cite{LeMiNe23} for recent progress concerning this question.

Another approach is to look at the Fokker-Planck equation associated to \eqref{eq:SDE0}. Given any test function $u_0$, $u(t,x)=\mathbb{E}(u_0(X_t)|X_0=x)$ is solution to the PDE
\be\label{eq:KFP}
\left\{
\begin{aligned}
    &\p_t u - \lll u = 0,\\
    &u_{|t=0}=u_0,
\end{aligned}
\right.
\ee
where
\bes
\lll = h\sum_{i,j}a_{i,j}\p_i\p_j + \sum_k\xi_k\p_k,
\ees
with $(a_{i,j})_{i,j} = \sigma\sigma^T$. Its adjoint problem is
\be\label{eq:KFPstar}
\left\{
\begin{aligned}
    &\p_t \psi - \lll^* \psi = 0,\\
    &\psi_{|t=0}=\psi_0,
\end{aligned}
\right.
\ee
where $A^*$ denotes the formal adjoint of any differential operator $A$. It is this last equation which is mainly called the Fokker-Planck equation whose solution is given by the density $\psi(t,X)$ that the random variable $X_t$ follows when it makes sense. Therefore studying $\lll$ or $\lll^*$ is a good way to obtain results concerning $X_t$.

A main question about this PDE is its resolution, implied by the maximal accretivity of $-\lll$ which is a real problem when it is not self-adjoint. In order to study the long time behaviour of \eqref{eq:KFP} an efficient strategy is to study the spectral properties of $\lll$ and particularly its smallest eigenvalues when they have non-negative real part. Determining its spectral gap informs us on equilibrium states and the metastability of \eqref{eq:SDE0}. When considering self-adjoint operators, the spectral Theorem directly implies the decreasing in time of the solution of \eqref{eq:KFP}, but it needs some more results otherwise. Although we will not go that far in this paper, let us mention that the non-self-adjoint setting yields consequent additional difficulties which are solved using the Gearhart-Prüss Theorem. In \cite{HeSj10_01} and \cite{HeSj21_01}, the authors have established a quantitative version of this Theorem uniform in $h$ which is the main argument to prove for example \cite[Corollaries 1.5 and 1.6]{BoLePMi22} or likewise \cite[Corollaries 1.8 and 1.9]{No23}.

The SDE \eqref{eq:SDE0} and its generator have been largely studied in the past decades, in particular when $h\to0$ and under more assumptions : for example taking $\xi = -\nabla V$ for $V:\R^d\to\R$ a potential and $\sigma=\Id$ we recover the overdamped Langevin process
\bes
dX_t = -\nabla V(X_t) + \sqrt{2h}dB_t,
\ees
whose generator is the Kramers-Smoluchowski operator 
\be\label{eq:KrSmo}
\lll = h\Delta - \nabla V\cdot\nabla
\ee
which is conjugated to the Witten Laplacian :
\bes
-e^{-V/2h}h\lll e^{V/2h} = \Delta_{\frac V2} := -h^2\Delta + \frac 14|\nabla V|^2 - \frac h2\Delta V.
\ees

Introduced in \cite{Wi82} in order to prove the Morse inequalities, this Laplacian is a non-negative self-adjoint operator that arises in many different domains such as control theory \cite{LaLe22} or dynamical system \cite{DaRi21}. As discussed previously, one main goal is to determine its spectrum (which we already know is included in $\R_+$) and especially the bottom of its spectrum for metastability questions and exit time estimates. At first, only estimates about the order of the bottom of the spectrum were proven in \cite{HeSj85_01}. They obtained that $\Delta_V$ has as much low lying eigenvalues as $V$ has minima and these eigenvalues satisfy the bound :
\bes
\lambda = O(e^{-c/h}),
\ees
for a certain $c>0$ uniform in $h$. In \cite{HoKuSt89} and \cite{Mi95}, the authors considered a smooth potential with small gradient on a compact manifold, they proved a law in-between Arrhenius and Eyring-Kramers. More precisely they obtained the following bounds on the eigenvalues
\be\label{eq:ArrEK}
a_jh^{5d} e^{-2S_j/h} \leq \lambda_j \leq b_jh^{-3d} e^{-2S_j/h}
\ee
for some $a_j,b_j > 0$ and explicit $S_j$. Progress on sharp asymptotics for these eigenvalues were slow due to topological restrains, but in \cite{BoEcGaKl04} and \cite{HeKlNi04_01} this barrier has been crossed and sharp estimates have been proven, obtaining the right order and describing precisely the prefactor
\bes
\lambda_j = z_jhe^{-2S_j/h}(1+o(1))
\ees
with explicit $z_j>0$ and $S_j>0$ and $1\leq j\leq n_0$ where $n_0$ denotes the number of minima of $V$. These kinds of formulae go back to pioneer work \cite{Kr40}, \cite{Er35} and are called Eyring-Kramers laws. Through functional analysis for self-adjoint operators, such estimates induce results about convergence rates of the semigroup associated to \eqref{eq:KrSmo}, the return to equilibrium being of the order $1/\lambda_2$ (because $\lambda_1=0$). More precisely, due to the different $S_j$, we seem to have a sort of stability of the system during exponentially large time intervals corresponding to the inverse of the $\lambda_j$, each one around the eigenfunction associated to $\lambda_j$, it is this phenomenon which is called metastability.

It has also been proven (see \cite{Da83}, \cite{FrWe84}, \cite{Ne21}, \cite{LeMiNe23}) that for this special operator, its first positive eigenvalue is the inverse of the mean exit time of the processus solving \eqref{eq:SDE0} in the boundary case.

These works and approaches which were doing just fine with this operator \eqref{eq:KrSmo} does not directly apply to non-self-adjoint operators which arise naturally from the homogeneous Langevin process. For example in $\R^{2d}$ taking $V:\R^d\to\R$ a potential, $\xi(x,v) = (v,-\nabla V(x)-v)$, $\gamma>0$ a friction coefficient and $\sigma = 0_x\oplus\Id_v$ we obtain
\be\label{eq:SDE1}
\left\{
\begin{aligned}
    dx_t &= v_tdt,\\
    dv_t &= (-\nabla V(x_t) - \gamma v_t)dt + \sqrt{2\gamma h}dB_t,
\end{aligned}
\right.
\ee
whose generator is the Kramers-Fokker-Planck operator
\bes
\lll = v\cdot\p_x - \p_xV\cdot\p_v + \gamma(h\Delta_v - v\cdot\p_v),
\ees
where $\Delta_v$ denotes the Laplacian acting only on the $v$ coordinates.

Due to the lack of self-adjointness because of the hamiltonian part $v\cdot\p_x-\p_xV\cdot\p_v$, the previous method had to be adapted, the main issue was resolvent estimates which were not free anymore. Through microlocal analysis, this problem was first solved in the non-semiclassical framework \cite{HeNi04} and then in the semiclassical one \cite{HeHiSj08-2}.
 
Equation \eqref{eq:SDE0} is used not only for modelling particle system, but its efficiency has been spread to other domains such as molecular dynamics \cite{VaBeMoCs13} or high dimensional data analysis \cite{ChChBaJo18}. The main advantage using the Langevin dynamic for numerical simulations instead of a usual Monte Carlo random walk is the use of the gradient of the potential $V$, which result in less wasted computation \cite{ScLeStCaCa06}, \cite{CaLeSt07}. Despite this benefit, it requires some precise information about that gradient which is a very challenging task. Thus the Adaptive Langevin dynamic was introduced in \cite{JoLe11}, \cite{LeSh16} from a fusion of a deterministic Nosé-Hoover scheme and a more usual overdamped Langevin process in order to reduce the needed knowledge about the gradient of $V$. By that time they mainly show their results through numerical simulations (we refer to \cite{LeRoSt10_01} for further details about the numerical study and modelisation of theses processes). The adaptive Langevin dynamic was next studied in \cite{LeSaSt} at fixed temperature where the authors determine some of its properties, namely spectral gap (see \cite[Theorem 2.1]{LeSaSt}) using hypocoercive estimates. The main question we want to address in this paper is to study how these properties depend on the semiclassical parameter $h$. In short, we model the gradient noise by another stochastic process which results in adding another unknown Brownian motion (that we can combine with the already existing one). But in order to retrieve the standard Gibbs state, we need to consider the friction coefficient to be a new variable, and all this leads to a slight modification in the SDE \eqref{eq:SDE1} :
\be\label{eq:SDE}
\left\{
\begin{aligned}
    dx_t &= v_tdt,\\
    dv_t &= (-\nabla V(x_t) - \nu y_tv_t-\gamma v_t)dt + \sqrt{2\gamma h}dB_t,\\
    dy_t &= \nu(|v_t|^2-dh)dt,
\end{aligned}
\right.
\ee
where $\gamma,\nu>0$ denotes positive parameters and the variable $X_t=(x_t,v_t,y_t)$ lives in $\R^d\times\R^d\times\R$. This SDE is \cite[(2.4)]{LeSaSt} with several names changed. Its associated generator called $\lll_{AdL}$ in \cite{LeSaSt} is
\be\label{eq:lll}
\lll = v\cdot\p_x-\p_xV\cdot\p_v+\nu((|v|^2-dh)\p_y-yv\cdot\p_v)+\gamma(h\Delta_v - v\cdot\p_v).
\ee

Even if it shares some similarities with the Langevin dynamics, \eqref{eq:SDE} and its generator \eqref{eq:lll} do not satisfy some of the crucial hypotheses made in \cite{BoLePMi22}. Mainly, they require that the operator is at most quadratic microlocally. In \eqref{eq:SDE} the terms that do not respect this assumption are $-\nu y_tv_tdt$ and $\nu|v_y|^2dt$ which will result in the cubic terms $\nu(|v|^2\p_y - yv\cdot\p_v)$ in \eqref{eq:lll}. That hypothesis inherited from \cite{HeHiSj08-2} was crucial for their microlocal estimates. In our work we manage to avoid that necessity by using the separated variable property of our particular operator which justifies the use of hypocoercivity methods in the spirit of \cite{DoMoSc15}. In that sense, \cite{BoLePMi22} is more general because it applies to a wide class of operator, but our work is not contained in theirs because of the operator's cubic term.

Therefore, this article is at the edge between \cite{LeSaSt} and \cite{BoLePMi22}, trying to use the arguments of the second reference in order to generalize the results of the first one in a semiclassical way and describe the low lying eigenvalues of their degenerate operator. Here we will obtain hypocoercivity, resolvent estimates and rough description of the eigenvalues uniform in the parameters $\gamma$ and $\nu$ (depending on $h$) but for sharp estimates, we had to fix $\gamma$ and $\nu$ independant of $h$.

\subsection{Statements.}
Considering $\lll$ from \eqref{eq:lll},we clearly have that $\lll1 =0$ and taking
\bes
f(x,v,y) = \frac{V(x)}{2}+\frac{|v|^2+y^2}{4},
\ees
one can show that  $\lll^* e^{-2f/h} = 0$. For this paper, we consider the conjugate operator $P = -e^{f/h}h\lll^* e^{-f/h}$. In our context, $\lll$ has only real coefficient, thus $\sigma(\lll^*) = \sigma(\lll)$, and we obtain
\bes
P =H_0+ \nu Y+\gamma \ooo,
\ees
where 
\bes
\left\{
\begin{aligned}
&H_0=v\cdot h\partial_x-\partial_xV\cdot h\partial_v,\\
&Y=(v h\partial_y-yh\partial_v)\circ v-h d(h\partial_y-\frac y 2),\\
&\ooo=-h^2\Delta_v+\frac{|v|^2}4-h\frac d 2.
\end{aligned}
\right.
\ees

We observe that we have the algebraic relations:
\be\label{eq:adj0}
H_0^*=-H_0,\;\;Y^*=-Y,\;\;\ooo^*=\ooo,
\ee
and inherited from the properties of $\lll$,
\be\label{eq:PPstarexp}
P(e^{-f/h}) = P^*(e^{-f/h}) = 0.
\ee

\begin{proposition}\label{prop:accretive}
The operator $P$ initially defined on $C^\infty_c(\R^{2d+1})$ admits a unique maximally accretive extension that we still denote by $(P,D(P))$, and we have $D(P) = \{u\in L^2(\R^{2d+1})\ |\ Pu\in L^2(\R^{2d+1})\}$.
\end{proposition}
We postpone the proof of this proposition to the Appendix.

\begin{assumption}
\label{ass.confin}
There exist $C>0$ and a compact  set $K\subset\R^d$ such that
$$V(x)\ \geq\  -C,\;\;\;
 \vert\nabla V(x)\vert \geq \frac 1 C\;\;\;\text{and}\;\;\; \vert\Hess V(x)\vert \leq C\,
$$
for all $x\in\R^d\setminus K$.
\end{assumption}
Under this assumption, it is known (see for example \cite[Lemma~3.14]{MeSc14}) that there exists $b>0$, such that $V(x)\geq  -b+b|x|$.

\begin{lemma}
Suppose that Assumption \ref{ass.confin} holds true. One has $e^{-f/h}\in D(P)$ and 
\be\label{eq:H_0YNmu}
H_0(e^{-f/h})=Y(e^{-f/h})=\ooo(e^{-f/h})=0.
\ee
\end{lemma}
\bp
The proof of \eqref{eq:H_0YNmu} is a simple computation, and therefore, we retrieve \eqref{eq:PPstarexp}, thus $e^{-f/h}\in D(P)$ thanks to the maximal accretivity of $P$ and $e^{-f/h}\in L^2(\R^{2d+1})$ thanks to Assumption \ref{ass.confin}.

\ep

\begin{assumption}\label{ass.morse}
The function $V$ is a Morse function.
\end{assumption}
Under Assumptions \ref{ass.confin} and \ref{ass.morse}, the set $\uuu$ of critical points of $V$ is finite. We denote by $\uuu^{(0)}$ the set of minima
of $V$ and $\uuu^{(1)}$ the set of critical points of index $1$. We shall also denote $n_0=\sharp\uuu^{(0)}$. As the critical points of $f$ are the $(x^*,0,0)$ for $x^*\in\uuu$, with the same index, we will identify those two and use $x^*$ instead of $(x^*,0,0)$ where it is clear which one we are really talking about ($x^*$ will mostly be denoted either $\m$ if of index $0$ or $\s$ if of index $1$).

Throughout the paper, we suppose that $V$ satisfies Assumptions \ref{ass.confin} and \ref{ass.morse}.

\begin{theorem}\label{thm}
There exists $h_0>0$, $c_0,c_1>0$, such that for all $h\in]0,h_0]$, there exists a subspace $G_h$ of $L^2(\R^{2d})$ of finite dimension $n_0 = \sharp\uuu^{(0)}$ such that for all $\gamma,\nu>0$, and all $u\in D(P)\cap G_h^\bot$, one has 
$$
\Vert (P-z)u\Vert_{L^2}\geq c_1g(h)\Vert u\Vert_{L^2}
$$
for every $z\in\C$ such that $\Re(z)\leq c_0g(h)$, where \be\label{eq:g}
g(h)=h\min\big(\nu^2h\gamma,\frac1\gamma,\frac{\gamma}{\nu^2h},\frac{\nu^2h}{\gamma}\big).
\ee
There exists an explicit constant $c_f>0$ depending only on $f$ such that if $g(h)$ satisfies
\be\label{eq:gthm2}
g(h)\geq e^{-\frac{\tilde c}{2h}}\mbox{ for any } \tilde c<c_f,
\ee
then there exists $\lambda_\m(h)\in\C$ for all $\m\in\uuu^{(0)}$ such that $\sigma(P)\cap\{\Re z\leq c_0g(h)\}=\{\lambda_\m(h),\m\in\uuu^{(0)}\}$ counted with multiplicity, and for all $\m\in\uuu^{(0)},\ |\lambda_\m(h)|\leq e^{-c_f/h}$. Moreover, for all $0<c_0'<c_1$
\bes
\forall |z|>c_0'g(h),\mbox{ such that }\Re z\leq c_0g(h),\ \vvert{(P-z)^{-1}}_{L^2}\leq\frac{2}{c_0'g(h)}.
\ees
\end{theorem}

\begin{remark}
Formally, when taking $h=1$, we recognize the conclusion of \cite[Corollary 1]{LeSaSt} noticing $\nu$ in our paper is $\eps^{-1}$ in theirs. Through similar hypocoercive methods we achieve to generalize their result to the semiclassical regime.
\end{remark}

This theorem, true in its general form will allow us to prove the following one, which describes a much more restrain case for the purpose of this paper : the double well. We will only consider this case because of its simplicity compared to the general one, the aim of this paper is to show that the sharp quasimodes and the methods developed in \cite{BoLePMi22} can be adapted to our operator although it does not satisfy some key assumption they made. To detail a bit more the technicality avoided here, in order to deal with a more generic case, one need to introduce several topological definitions regarding the minima of $V$ and sets around theses minima that will be essential to defined sharp quasimodes and have the most precise estimates. Moreover, at the end we obtain a matrix whose eigenvalues are the eigenvalues we are looking for, but in our case it is a mere $2\times 2$ matrix with three zeros. In the general case, the matrix is not even diagonal and it needs a non trivial study to extract its eigenvalues. We will get into the general case in a forecasting paper.

\begin{theorem}\label{thm3}
Let us suppose $\uuu^{(0)}=\{\underline{\m},\widehat{\m}\}$ where $\underline{\m}$ is the unique global minimum of $V$, $\uuu^{(1)}=\{\s\}$ and $\gamma,\nu >0$ are fixed. There exists $c_0,h_0>0$, such that for all $h\in]0,h_0]$, one has
\bes
\sigma(P)\cap\{\Re z\leq c_0g(h)\} = \{0,\lambda\},
\ees
where $g(h)\propto h^2$ is as in \eqref{eq:g} and with
\bes
\lambda=\frac{\mu(\s)(\det\Hess V(\widehat{\m}))^{\frac12}}{2\pi|\det\Hess V(\s)|^{\frac12}}he^{-S(\widehat{\m})/h}(1+O(\sqrt h)),
\ees
where $\mu(\s) = \frac{1}{2}(-\gamma+\sqrt{\gamma^2+4\eta})>0$ with $-\eta$ the sole negative eigenvalue of $\Hess_\s V$ and
\bes
S(\widehat{\m})=V(\s)-V(\widehat{\m}).
\ees
\end{theorem}

\begin{remark}
If familiar with this sort of results, one should expect to have a factor $2$ befront the $S$, but since we take heights in terms of $V$ and not $f$, we do not have that $2$.
\end{remark}

\begin{remark}
One can prove a similar theorem without the double well assumption, it requires much more geometric constructions. Going further in the development of $w$ (defined in \eqref{eq:w}) in order to obtain higher principal orders in \eqref{eq:Pchi}, one can transform the $1+O(\sqrt h)$ into $1+O(h)$ and even obtain a full semiclassical asymptotic development $1+h\sum_{j\geq0}a_jh^j$ with explicit $a_j$. Consult \cite{BoLePMi22} for more details and an explicit way to obtain that generalisation.
\end{remark}

In the next section, we will show some hypocoercive estimates for $P$ following the work of \cite{DoMoSc15} by defining an adapted auxiliary operator leading to the proof of Theorem \ref{thm} giving rough localization on the spectrum of $P$ and a resolvent estimate. In order to obtain precise estimate of the eigenvalues of $P$ we will first make kinds of WKB constructions in Section \ref{sec:LocalQuasimode} to resolve $Pu=0$ locally around a given saddle point following \cite{BoLePMi22}'s method. These constructions will help us define good quasimodes globally in Section \ref{sec:TwoWells} for a rather simple example of potential, which will be key to the proof of Theorem \ref{thm3}.

\subsection*{Acknowledgements.} The author is grateful to Laurent Michel for his advice through this work and to Gabriel Stoltz for helpful discussions.\newline This work is supported by the ANR project QuAMProcs 19-CE40-0010-01.

\section{Hypocoercive estimates}

The goal of this section is to prove some hypocoercive estimates that will lead to Theorem \ref{thm}. This is achieved by comparing $P$ to a well-known operator, the Witten Laplacian. We follow \cite{DoMoSc15}'s method, and therefore we define an auxiliary operator, function of the skew-adjoint part of $P$ and a projector onto the kernel of the self-adjoint part of $P$. This way it shall "contain" all the information about $P$ and we will able to use it to prove the hypocoercivity.

\subsection{Witten Laplacians and an auxiliary operator.}
We thus first consider the semiclassical Witten Laplacian associated to the function $\frac V2$, acting on $L^2(\R^d_x)$
\bes
\Delta_{\frac V2}=-h^2\Delta_x+\frac 14 |\nabla V|^2-\frac h 2\Delta V,
\ees
and the semiclassical Witten Laplacian associated to the function $y\mapsto \frac{y^2}4$, acting on $L^2(\R_y)$
\bes
N_y=-h^2\partial_y^2+\frac{y^2}4-\frac h 2.
\ees
Throughout the paper, we denote $\delta_{x_i}=h\partial_{x_i}+\frac{\partial_{x_i}V}2$ and $\delta_{y}=h\partial_{y}+\frac{y}2$ the associated twisted derivatives. One has the identities
$$
\Delta_{\frac V2}=\sum_{i=1}^d\delta_{x_i}^*\delta_{x_i}\text{ and }N_y=\delta_{y}^*\delta_{y}.
$$
Along with $\delta_v = h\p_v + \frac{v}{2}$ which then gives $\ooo = \delta_v^*\delta_v$, these twisted derivatives allow us to rewrite $H_0$ and $Y$ in a more fancy way through direct computations
\be\label{eq:H_0Ydelta}
\begin{split}
    &H_0 = v\cdot\delta_x - \p_xV\cdot\delta_v,\\
    &Y = (|v|^2-dh)\delta_y - yv\cdot\delta_v.
\end{split}
\ee

We now introduce the function $\rho(v)=(2\pi h)^{-\frac d 4}e^{-\frac{|v|^2}{4h}}$ and the orthogonal projector onto the kernel of $\ooo$ defined on $L^2(\R^{2d+1})$ by
$$
\Pi_\rho u(x,v,y)=\int_{\R^d}u(x,v',y)\rho(v')dv' \rho(v)=u_\rho(x,y)\rho(v),
$$
where we denoted
\bes
u_\rho = \< u, \rho\>_{L^2_v(\R^d)}.
\ees

Let us denote  $Z=H_0+ \nu Y$ the skew-adjoint part of $P$ and notice that we have
\bes
\Pr Z\Pr = \ooo\Pr = 0.
\ees
Indeed using \eqref{eq:H_0YPi}, for any $u\in L^2(\R^{2d+1})$,
\bes
\begin{split}
    \Pr Z\Pr &= \Pr(v\cdot\delta_x+\nu(|v|^2-dh)\delta_y)\Pr\\
    \Pr v_j\Pr u &= cu_\rho\rho\int v_je^{-\frac{|v|^2}{2h}}\dv = 0\\
    \Pr (|v|^2-dh)\Pr u &= chu_\rho\rho\int (|v|^2-d)e^{-\frac{|v|^2}{2}}\dv = 0
\end{split}
\ees
whence $\Pr Z\Pr = 0$ (with $c$ a constant that changed from line two to three). Moreover we have the following lemma that will be useful many times in the following.

\begin{lemma}\label{lem:bornVPi}
For any $j=1,\ldots, d$, the operator $v_j\Pi_\rho$ is bounded on $L^2$ and 
$$
\forall k\in\N,\ \ \Vert v_j^k\Pi_\rho\Vert_{L^2\rightarrow L^2}=O(h^{k/2}).
$$
\end{lemma}

\bp
We notice that for $u\in L^2(\R^{2d+1})$, $v_j\Pr u = -2h\p_{v_j}\Pr u$ hence the result.

\ep

We define for $\alpha>0$
\be\label{eq:defA}
A=\big(h\alpha+h^{-1}(Z\Pi_\rho)^*(Z\Pi_\rho)\big)^{-1}(Z\Pi_\rho)^*.
\ee
This auxiliary operator is introduced in \cite{DoMoSc15} and used in \cite{LeSaSt} in order to ease the calculus in the proof of Theorem \ref{thm}. This kind of method to compute hypocoercivity was mainly introduced and used at first in \cite{Vi09}, \cite{HeNi04} and \cite{He06}.
\begin{lemma}\label{lem:Aborne}
The operator $A$ is bounded on $L^2(\R^{2d+1})$, it satisfies 
$$A = \Pr A = A(1-\Pr)$$  and one has the estimate
$$
\Vert A\Vert_{L^2}\leq \frac{1}{\sqrt{\alpha}}.
$$
\end{lemma}
\bp
The equation $Au = w$ is equivalent to
\bes
(Z\Pr)^*u = h\alpha w + h^{-1}(Z\Pi_\rho)^*(Z\Pi_\rho) w.
\ees
Writing this as $h\alpha w = h^{-1}\Pr Z^2\Pr w - \Pr Z u$ proves $A = \Pr A$. And the equality $A=A(1-\Pr)$ comes from $\Pr Z\Pr=0$ which has been proven before stating the lemma.

For the bound, we have
\bes
AA^* = \big(h\alpha+h^{-1}(Z\Pi_\rho)^*(Z\Pi_\rho)\big)^{-1}(Z\Pi_\rho)^*(Z\Pi_\rho)\big(h\alpha+h^{-1}(Z\Pi_\rho)^*(Z\Pi_\rho)\big)^{-1},
\ees
and through functional calculus we know that
\bes
\vvert{\big(h\alpha+h^{-1}(Z\Pi_\rho)^*(Z\Pi_\rho)\big)^{-1}(Z\Pi_\rho)^*(Z\Pi_\rho)} \leq \sup_{x\geq0} \frac{x}{h\alpha + h^{-1}x} = h
\ees
and
\bes
\vvert{\big(h\alpha+h^{-1}(Z\Pi_\rho)^*(Z\Pi_\rho)\big)^{-1}} \leq \sup_{x\geq0} \frac{1}{h\alpha + h^{-1}x} = (h\alpha)^{-1},
\ees
hence the bound.

\ep

And because $\Pr$ is a projector on the kernel of $\delta_v$, we thus have
\be\label{eq:H_0YPi}
H_j\Pi_\rho=v_j\delta_{x_j}\Pi_\rho\;\;\text{ and } \;\; Y\Pi_\rho=(|v|^2-dh)\delta_y\Pi_\rho,
\ee
where $H_j = v_jh\p_{x_j}-\p_{x_j}Vh\p_{v_j}$ and $H_0 = \sum_{j=1}^dH_j$. We also recall the commutation rules
\be\label{eq:comutWit0}
[\delta_{x_i},\delta_{x_j}]=0,\;\;[\delta_{x_i},\delta_{x_j}^*]=h\partial^2_{ij} V,\;\;[\delta_{x_i}, \delta_{x_j}^*\delta_{x_j}]= h\partial^2_{ij}V\delta_{x_j}
\ee
and 
\be\label{eq:comutWit2}
[\delta_{y},N_y]=h \delta_{y}.
\ee

We can now state a link between our auxiliary operator $A$ and a Witten Laplacian.

\begin{lemma}\label{lem:compsquare}
One has
\be\label{eq:lemBPirho}
(Z\Pi_\rho)^*(Z\Pi_\rho)=hB\Pi_\rho,
\ee
where
$$
B=\Delta_{\frac V 2}+2d\nu^2hN_y.
$$
\end{lemma}
\bp
We have $Z\Pr = (H_0+\nu Y)\Pr = (v\cdot\delta_x+\nu(|v|^2-dh)\delta_y)\Pr$ but, because of the parity of $\rho$, for all $u,w\in\ccc^\infty_c(\R^{2d+1})$,
\bes
\begin{split}
    \< Z\Pr u, Z\Pr w\> &= \< v\cdot\delta_x\Pr u,v\cdot\delta_x\Pr w\> + \nu^2\< (|v|^2-dh)\delta_y\Pr u,(|v|^2-dh)\delta_y\Pr w\>\\&\hspace{1.4cm}+2\nu\Re\underbrace{\< v\cdot\delta_x\Pr u,(|v|^2-dh)\delta_y\Pr w\>}_{=0},\\
\end{split}
\ees
the last scalar product involves an integral over $\R^d$ of an odd function of $v$ it is therefore null. With the same argument, in the double sum, we only have the diagonal terms :
\bes
\begin{split}
    \< v\cdot\delta_x\Pr u,v\cdot\delta_x\Pr  w\> &= \sum_{i=1}^d \< v_i\delta_{x_i}\Pr u,v_i\delta_{x_i}\Pr  w\> = \sum_{i=1}^d \< v_i^2\delta_{x_i}^*\delta_{x_i}\Pr u,\Pr w\>\\
    &= \sum_{i=1}^d \< v_i^2\delta_{x_i}^*\delta_{x_i}u_\rho\rho,w_\rho\rho \> = \sum_{i=1}^d \<\delta_{x_i}^*\delta_{x_i}u_\rho,w_\rho \>_{L^2_{x,y}} \< v_i^2\rho,\rho \>_{L^2_v}\\
    &= \sum_{i=1}^d \<\delta_{x_i}^*\delta_{x_i}u_\rho\rho,w \>_{L^2_{x,v,y}} \< v_i^2\rho,\rho \>_{L^2_v}.
\end{split}
\ees
By integration by parts, we note
\bes
\< v_i^2\rho,\rho \>=(2\pi)^{-\frac d2}\int_{\R^d}v_i^2e^{-\frac{|v|^2}{2h}}\frac{dv}{h^\frac{d}{2}}=h
\ees
and so
\bes
\< v\cdot\delta_x\Pr u,v\cdot\delta_x\Pr w \> = \< h\Delta_{\frac{V}{2}}\Pr u, w \>.
\ees
With very similar computations, we obtain
\bes
\begin{split}
    \< (|v|^2-dh)\delta_y\Pr u,(|v|^2-dh)\delta_y\Pr w \> &= \< \delta_y^*\delta_y u_\rho\rho,w \> \< (|v|^2-dh)^2\rho,\rho \>\\
    \< (|v|^2-dh)\delta_y\Pr u,(|v|^2-dh)\delta_y\Pr w \> &= \< dh^2N_y\Pr u, w \>,
\end{split}
\ees
noticing that $\< (|v|^2-dh)^2\rho,\rho \> = dh^2$.

Finally, we observe that we indeed have proved \eqref{eq:lemBPirho}
\bes
(Z\Pi_\rho)^*(Z\Pi_\rho)=hB\Pi_\rho.
\ees

\ep

One direct consequence of Lemma \ref{lem:compsquare} and \eqref{eq:defA} is that we have
\bes
A=(h\alpha+B\Pr)^{-1}(Z\Pr)^*.
\ees

We can now use the well-known properties of the Witten Laplacian to obtain a lower bound on $B\Pr$ on the orthogonal of a finite dimensional space. Let $\chi_\m$, $\m\in \uuu^{(0)}$ be some cutoffs in $\ccc_c^\infty(\R^d)$ such that $\chi_\m$ is supported in $B(\m,r)$ for some $r>0$ to be chosen small enough and $\chi_\m$ is constant near $\m$. We then introduce the quasimodes
\bes
f_\m(x,v,y)=\chi_\m(x) e^{-(f(x,v,y)-f(\m))/h},
\ees
and we set the constant $\chi_\m(\m)$ such that $f_\m$ is of norm one in $L^2(\R^{2d+1})$. Thus, with a Laplace method we observe that $\chi_\m(\m)$ is of order $h^{-\frac d2-\frac14}$.

For $r>0$ small enough, these functions have disjoint support and hence the vector space 
\bes
F_h=\span\{f_\m,\;\m\in\uuu^{(0)}\}
\ees
has dimension $n_0$. We in fact have that $G_h$ in Theorem \ref{thm} is $F_h$ we just defined. It is a natural space to consider noticing that $e^{-f/h}\R$ is the kernel of $B\Pr$ which should not be surprising since $B\Pr$ is a self-adjoint operator built to behave like $P$.

\begin{lemma}\label{eq:minorB}
There exists $c_0,h_0>0$ such that for all $h\in]0,h_0]$, $\nu>0$ and $u\in F_h^\bot$,  one has
\bes
\<B\Pi_\rho u,u\>\geq c_0h\min(1,\nu^2h)\Vert \Pi_\rho u\Vert^2.
\ees
\end{lemma}

\bp
We set $W(x,y)=\frac{V(x)}{2}+\frac{y^2}{4}$, hence $\Delta_W=\Delta_{\frac{V}{2}}+N_y$, and we see that $W$ has the same property as $V$ : if $V\geq-C$ then so is $W$, $|\nabla W|^2=\frac{1}{4}(|\nabla V|^2+y^2)$ and $\Hess_W(x,y)=\frac{1}{2}\left(\begin{matrix}
    \Hess_V(x)&0\\
    0&1\\
\end{matrix}\right)$. Therefore $W$ satisfies Assumptions \ref{ass.confin} and \ref{ass.morse} as much as $V$, and the minima of $W$ are the $(\m,0)$ where $\m\in\uuu^{(0)}$. In order to lighten the notations we will identify $\m$ and $(\m,0)$, likewise we will identify $\uuu^{(0)}$ with $\uuu^{(0)}\times\{0\}$. We also denote $\delta_W=h\nabla+\nabla W$, and $X=(x,y)\in\R^{d+1}$.

Using known facts about the Witten laplacian (see for example \cite[Theorem 11.1]{CyFrKiSi87_01} or \cite{HeSj85_01} for the exponential bound) we have that
\be\label{eq:gap}
\exists c,\eps,h_0>0,\forall h\in\,]0,h_0]\ \sigma(\Delta_W)\cap\,]e^{-c/h},\eps h[\,=\emptyset,
\ee
and $\Delta_W$ has exactly $n_0$ eigenvalues in $[0,e^{-c/h}]$ that we denote $E_n(\Delta_W)$, where $1\leq n \leq n_0$.

We also denote $\widetilde{F}_h=\{u_\rho, u\in F_h\}=\span(\widetilde{f}_\m)_{\m\in\uuu^{(0)}}$ where $\widetilde{f}_\m = f_\m(\cdot,0,\cdot)$, then we will admit for now that because of \eqref{eq:gap},
\bes
\exists\eps'>0,\ \forall u\in\widetilde{F}_h^\bot,\ \<\Delta_Wu,u\>\geq\eps' h\vvert{u}^2.
\ees
Therefore, for $u\in F_h^\bot$
\bes
\begin{split}
    \< B\Pr u,u\> &= \< B\Pr^2 u,u\> = \<\Pr B\Pr u,u\> = \< B\Pr u,\Pr u\>\\
    &= \<\Delta_{\frac{V}{2}}\Pr u,\Pr u\> + 2d\nu^2h\< N_y\Pr u,\Pr u\>\\
    &\geq \min(1,2d\nu^2h)\<\Delta_W\Pr u,\Pr u\> \geq \eps'h\min(1,\nu^2h)\vvert{\Pr u}^2.
\end{split}
\ees

\ep

And so we proved the lemma, let us now show what we have admitted : 
\begin{lemma}\label{lem:hypocoerW}
There exists $\eps'>0,$ such that for all $u\in\widetilde{F}_h^\bot,\ \<\Delta_Wu,u\>\geq\eps' h\vvert{u}^2.$
\end{lemma}

\bp We first define the Riesz projector on the eigenvectors associated to the small eigenvalues : $\Pi_W=\frac{1}{2i\pi}\int_{\p\negthinspace D}(z-\Delta_W)^{-1}\dz$ where we denote $\p\negthinspace D = \p\negthinspace D(0,\frac\eps2 h)$ the circle centered in $0$ of radius $\frac\eps2 h$ positively oriented, where $\eps$ is defined in \eqref{eq:gap}. Thus,
\bes
\Pi_W-\Id=\frac{1}{2i\pi}\int_{\p\negthinspace D}((z-\Delta_W)^{-1}-z^{-1})\dz=\frac{1}{2i\pi}\int_{\p\negthinspace D}(z-\Delta_W)^{-1}\Delta_Wz^{-1}\dz,
\ees
applied to the $\widetilde{f}_\m$, we get
\be\label{eq:PiWf}
\Pi_W \widetilde{f}_\m-\widetilde{f}_\m=\frac{1}{2i\pi}\int_{\p\negthinspace D}\underbrace{(z-\Delta_W)^{-1}}_{=O(h^{-1})}\underbrace{\Delta_W(\widetilde{f}_\m)}_{=O(e^{-c/h})}\frac\dz z=O(e^{-c'/h}).
\ee

Then by the spectral theorem (which we can use because $\Delta_W$ is self-adjoint), noting $\phii_n$ normalized eigenvectors of $\Delta_W$ associated to $E_n(\Delta_W)$, we have for $u\in D(\Delta_W)$
\bes
\begin{split}
     \< \Delta_W u,u\> &= \sum_{n\leq n_0}E_n(\Delta_W)|\< u,\phii_n\>|^2+\int_{\eps h}^\infty\lambda d\<E_\lambda u,u\>\\
     &\geq \int_{\eps h}^\infty\lambda d\<E_\lambda u,u\> \geq \eps h\bigg(\vvert{u}^2-\sum_{n\leq n_0}|\< u,\phii_n\>|^2\bigg).
\end{split}
\ees
We now want to show that $\exists c>0,\forall u\in\widetilde{F}_h^\bot,\ \vvert{u}^2-\sum_{n\leq n_0}|\langle u,\phii_n\rangle|^2\geq c\vvert{u}^2$, or in an equivalent way $\sum_{n\leq n_0}|\langle u,\phii_n\rangle|^2\leq c'\vvert{u}^2$ with $c'<1$. For $\m\in\uuu^{(0)}$, we have that $\Vert\widetilde{f}_\m\Vert_{L^2(\R^{d+1})}$ is of order $h^{-\frac d4}$ because $\Vert f_\m\Vert_{L^2(\R^{2d+1})} = 1$. Moreover, because of \eqref{eq:PiWf}, we have that $\span(\Pi_W\widetilde{f}_\m)_{\m\in\uuu^{(0)}} = \span(\phii_n)_{n\leq n_0}$, therefore there exists $a_{n,\m}\in\R$ of order $h^{\frac d4}$ such that for all $n\leq n_0$, $\phii_n = \sum_{\m\in\uuu^{(0)}}a_{n,\m}\Pi_W\widetilde{f}_\m$. Using \eqref{eq:PiWf} we obtain
\bes
\begin{split}
    \sum_{n\leq n_0}\<\cdot,\phii_n\>\phii_n &= \sum_{n\leq n_0}\sum_{\m,\m'\in\uuu^{(0)}}a_{n,\m}a_{n,\m'}\<\cdot,\Pi_W\widetilde f_\m\>\Pi_W\widetilde f_{\m'}\\
    &= \sum_{n\leq n_0}\sum_{\m,\m'\in\uuu^{(0)}}a_{n,\m}a_{n,\m'}\<\cdot,\widetilde f_\m\>\widetilde f_{\m'} + O(e^{-c/h}).
\end{split}
\ees
Now considering that $\forall u \in \widetilde{F}_h^\bot,\forall \m\in\uuu^{(0)}$, $\<u,\widetilde{f}_\m\> = 0$, we have that there exists $C>0$ such that for all $h>0$ small enough and all $u\in\widetilde{F}_h^\bot$,
\bes
\sum_{n\leq n_0}|\langle u,\phii_n\rangle|^2\leq e^{-C/h}\vvert{u}^2.
\ees

\ep

\subsection{Boundedness of remaining terms.}
Now, in order to apply \cite{DoMoSc15}'s method, we need to bound some remaining terms that appear, this will be the role of Lemma \ref{lem:majorJprime}.

\begin{lemma}\label{lem:computprod}
One has the following identities
\be\label{prod1}
H_i^*H_j\Pi_\rho=(-v_iv_j\delta_{x_i}\delta_{x_j}+\delta_{i,j} h\partial_{x_i} V\delta_{x_j})\Pi_\rho.
\ee
\be\label{prod2}
Y^*Y\Pi_\rho=\big((|v|^2-dh)^2\delta_y^*\delta_y-((|v|^2-dh)^2-2h|v|^2)y\delta_y\big)\Pi_\rho.
\ee
\be\label{prod3}
H_i^*Y\Pi_\rho=\big(-v_i(|v|^2-dh)\delta_{x_i}\delta_y+2v_ih\partial_{x_i}V \delta_y\big)\Pi_\rho.
\ee
\be\label{prod4}
Y^*H_i\Pi_\rho=\big(-v_i(|v|^2-dh)\delta_{x_i}\delta_y+v_ihy \delta_{x_i}\big)\Pi_\rho.
\ee
\end{lemma}

These identities will be useful for the following lemma but since its proof is mere calculus we postpone it to the Appendix.

\begin{lemma}\label{lem:majorJprime}
There exists $C,h_0>0$ such that for all $h\in]0,h_0]$,  $\nu>0$ and for all $u\in\ccc^\infty_c(\R^{2d+1})$, one has
\be\label{eq:majJ1}
|\<AZ (1-\Pi_\rho)u,u\>|\leq C(1+\nu\sqrt{h}\alpha^{-\frac12}+\alpha^{-\frac12})h\Vert \Pi_\rho u\Vert\, \Vert(1- \Pi_\rho) u\Vert.
\ee
\be\label{eq:majJ2}
|\<A\ooo u,u\>|\leq C\alpha^{-\frac12}h\Vert \Pi_\rho u\Vert\, \Vert(1- \Pi_\rho) u\Vert.
\ee
\be\label{eq:majJ3}
|\<Z u,A u\>|\leq Ch\Vert(1- \Pi_\rho) u\Vert^2.
\ee
\end{lemma}
\bp
Within this proof, $C$ will denote a positive constant that may only depends on the dimension $d$ and may change from line to line, and $u$ is a test function in $\ccc^\infty_c(\R^{2d+1})$.

Let us start with the proof of \eqref{eq:majJ1}. Since $A=\Pr A$, by the Cauchy-Schwarz inequality, we have
\bes
|\<AZ (1-\Pi_\rho)u,u\>| = |\<\Pr AZ (1-\Pi_\rho)u,u\>| \leq \vvert{AZ}\vvert{\Pr u}\vvert{(1-\Pr)u}.
\ees
Therefore, bounding the operator $AZ$ (or equivalently its adjoint) by $C(1+\nu\sqrt{h}\alpha^{-\frac12}+\alpha^{-\frac12})h$ is enough to prove \eqref{eq:majJ1}. One has
\bes
\begin{aligned}
Z^*A^* &= Z^*Z\Pi_\rho(h\alpha+B\Pr)^{-1}\\
&= \Big(\sum_{i,j}H_i^*H_j+\nu^2Y^*Y+\nu\sum_i(Y^*H_i+H_i^*Y)\Big)(h\alpha+B\Pr)^{-1}\Pi_\rho
\end{aligned}
\ees
and we estimate each term separately.
We start with the term involving $H_i^*H_j$. From \eqref{prod1}, we deduce that 
\be\label{mJp1}
H_i^*H_j(h\alpha+B\Pr)^{-1}\Pr = (-v_iv_j\delta_{x_i}\delta_{x_j}+\delta_{i,j} h\partial_{x_i} V\delta_{x_j})(h\alpha+B\Pr)^{-1}\Pi_\rho.
\ee
First we observe that
\bes
\begin{split}
\Vert\delta_{x_i}\delta_{x_j}(h\alpha+\Delta_{\frac{V}{2}})^{-1}u\Vert^2 &\leq \sum_{k,l}\Vert\delta_{x_k}\delta_{x_l}(h\alpha+\Delta_{\frac{V}{2}})^{-1}u\Vert^2\\
&\leq \sum_{k,l}\<\delta_{x_k}^*\delta_{x_k}\delta_{x_l}(h\alpha+\Delta_{\frac{V}{2}})^{-1}u,\delta_{x_l}(h\alpha+\Delta_{\frac{V}{2}})^{-1}u\>\\
&\leq \sum_l\<\Delta_{\frac{V}{2}}\delta_{x_l}(h\alpha+\Delta_{\frac{V}{2}})^{-1}u,\delta_{x_l}(h\alpha+\Delta_{\frac{V}{2}})^{-1}u\>\\
&\leq  \sum_l\<\Delta_{\frac{V}{2}}(h\alpha+\Delta_{\frac{V}{2}})^{-1}u,\delta_{x_l}^*\delta_{x_l}(h\alpha+\Delta_{\frac{V}{2}})^{-1}u\>\\
&\phantom{*****} +\sum_l\<[\Delta_{\frac{V}{2}},\delta_{x_l}](h\alpha+\Delta_{\frac{V}{2}})^{-1}u,\delta_{x_l}(h\alpha+\Delta_{\frac{V}{2}})^{-1}u\>\\
&\leq \Vert \Delta_{\frac{V}{2}}(h\alpha+\Delta_{\frac{V}{2}})^{-1}u\Vert^2\\
&\phantom{*****} -\sum_{k,l}\<h\partial^2_{kl}V\delta_{x_k}(h\alpha+\Delta_{\frac{V}{2}})^{-1}u,\delta_{x_l}(h\alpha+\Delta_{\frac{V}{2}})^{-1}u\>\\
&\leq C(1+h\max_k\Vert \delta_{x_k}(h\alpha+\Delta_{\frac{V}{2}})^{-\frac 12}\Vert^2\Vert (h\alpha+\Delta_{\frac{V}{2}})^{-\frac 12}\Vert^2)\Vert u\Vert^2,
\end{split}
\ees

where we used \eqref{eq:comutWit0} to compute the commutator, Assumption \ref{ass.confin} and Lemma \ref{lem:resolvestim} to get the last estimates. Thanks to  Lemma \ref{lem:elliptreg} and Lemma \ref{lem:resolvestim}, this implies 
\be\label{mJp2}
\Vert\delta_{x_i}\delta_{x_j}(h\alpha+\Delta_{\frac{V}{2}})^{-1}u\Vert \leq C (1+h^\frac12h^{-\frac12}\alpha^{-\frac12}) \Vert u\Vert \leq C(1+\alpha^{-\frac12}) \Vert u\Vert.
\ee
Using Lemma \ref{lem:bornVPi} and Lemma \ref{lem:bornres} this implies
\be\label{mJp2bis}
v_iv_j\delta_{x_i}\delta_{x_j}(h\alpha+B\Pr)^{-1}\Pi_\rho = (1+\alpha^{-\frac12})O(h).
\ee

Similarly, since $|\nabla V|^2\leq 4(\Delta_{\frac{V}{2}}+\frac h2\Delta V)$ in the sense of operators, we have
\bes
\begin{split}
\Vert\partial_{x_i} V\delta_{x_j}(h\alpha+d\Delta_{\frac{V}{2}})^{-1}u\Vert^2 &\leq \sum_{k,l}\<|\partial_{x_k} V|^2\delta_{x_l}(h\alpha+d\Delta_{\frac{V}{2}})^{-1}u,\delta_{x_l}(h\alpha+d\Delta_{\frac{V}{2}})^{-1}u\>\\
&\leq C\sum_l\<\Delta_{\frac{V}{2}}\delta_{x_l}(h\alpha+d\Delta_{\frac{V}{2}})^{-1}u,\delta_{x_l}(h\alpha+d\Delta_{\frac{V}{2}})^{-1}u\>\\
&\phantom{********}+Ch\max_k\Vert |\Delta V|^{\frac 12} \delta_{x_k}(h\alpha+d\Delta_{\frac{V}{2}})^{-1}u\Vert^2\\
&\leq C\sum_l\<\Delta_{\frac{V}{2}}\delta_{x_l}(h\alpha+d\Delta_{\frac{V}{2}})^{-1}u,\delta_{x_l}(h\alpha+d\Delta_{\frac{V}{2}})^{-1}u\>\\
&\phantom{********}+Ch\max_k\Vert  \delta_{x_k}(h\alpha+d\Delta_{\frac{V}{2}})^{-1}u\Vert^2,
\end{split}
\ees
using Assumption \ref{ass.confin}, which implies by the same arguments as above that 
\bes
\Vert\partial_{x_i} V\delta_{x_j}(h\alpha+\Delta_{\frac{V}{2}})^{-1}u\Vert\leq C(1+\alpha^{-\frac12})\Vert u\Vert.
\ees
Using again  Lemma \ref{lem:bornres}, it follows that
\be\label{mJp3}
h\partial_{x_i} V\delta_{x_j}(h\alpha+B\Pr)^{-1}\Pi_\rho = (1+\alpha^{-\frac12})O(h).
\ee
Combining \eqref{mJp1}, \eqref{mJp2bis} and \eqref{mJp3} , we finally get
\be\label{eq:finalH_0H_0}
H_i^*H_j(h\alpha+B\Pr)^{-1}\Pr = (1+\alpha^{-\frac12})O(h).
\ee

Now, from \eqref{prod2}, we have
\bes
\begin{split}
    \nu^2Y^*Y(h\alpha+B\Pr)^{-1}\Pr &= \nu^2\big((|v|^2-dh)^2\delta_y^*\delta_y-((|v|^2-dh)^2-2h|v|^2)y\delta_y\big)\\
    &\phantom{******}\times(h\alpha+B\Pr)^{-1}\Pr.
\end{split}
\ees
First we notice that due to Lemma \ref{lem:resolvestim},
\bes
\vvert{\delta_y^*\delta_y(h\alpha+2d\nu^2hN_y)^{-1}} = \vvert{N_y(h\alpha+2d\nu^2hN_y)^{-1}}\leq C\nu^{-2}h^{-1}.
\ees
And using Lemma \ref{lem:bornVPi} and \ref{lem:bornres}, we then get :
\bes
\begin{split}
    \vvert{(|v|^2-dh)^2\delta_y^*\delta_y(h\alpha+B\Pr)^{-1}\Pr} &\leq \vvert{(|v|^2-dh)^2\Pr}\vvert{\delta_y^*\delta_y(h\alpha+2d\nu^2hN_y)^{-1}}\\
    &\phantom{******} \times\vvert{(h\alpha+2d\nu^2hN_y)(h\alpha+B\Pr)^{-1}}\\
    &\leq C\nu^{-2}h.
\end{split}
\ees
Hence
\be\label{mJp8}
(|v|^2-dh)^2\delta_y^*\delta_y(h\alpha+B\Pr)^{-1}\Pr = \nu^{-2}O(h).
\ee

Now since $y^2\leq 4N_y+2h$ in the sense of operators, we have :
\bes
\begin{split}
    (\ast_1) &:= \vvert{y\delta_y(h\alpha+2d\nu^2hN_y)^{-1}u}^2\\
    &= \< y^2\delta_y(h\alpha+2d\nu^2hN_y)^{-1}u,\delta_y(h\alpha+2d\nu^2hN_y)^{-1}u\>\\
    &\leq 4\< N_y\delta_y(h\alpha+2d\nu^2hN_y)^{-1}u,\delta_y(h\alpha+2d\nu^2hN_y)^{-1}u\>\\
    &\phantom{*******} +2h\vvert{\delta_y(h\alpha+2d\nu^2hN_y)^{-1}u}^2\\
    &= 4\vvert{N_y(h\alpha+2d\nu^2hN_y)^{-1}u}^2+2h\vvert{\delta_y(h\alpha+2d\nu^2hN_y)^{-1}u}^2\\
    &\phantom{*******}+4\<[N_y,\delta_y](h\alpha+2d\nu^2hN_y)^{-1}u,\delta_y(h\alpha+2d\nu^2hN_y)^{-1}u\>\\
    &= 4\vvert{N_y(h\alpha+2d\nu^2hN_y)^{-1}u}^2-2h\vvert{\delta_y(h\alpha+2d\nu^2hN_y)^{-1}u}^2\\
    &\leq 4\nu^{-4}h^{-2}\vvert{u}^2,
\end{split}
\ees
where we used \eqref{eq:comutWit2} to compute the commutator and Lemma \ref{lem:elliptreg} for the last estimate. Thus, thanks to Lemma \ref{lem:bornVPi} and \ref{lem:bornres}, it implies :
\bes
\begin{split}
    (\ast_2) &:= \vvert{((|v|^2-dh)^2-2h|v|^2)y\delta_y(h\alpha+B\Pr)^{-1}\Pr}\\ &\leq \vvert{y\delta_y(h\alpha+2d\nu^2hN_y)^{-1}}\vvert{(h\alpha+2d\nu^2hN_y)(h\alpha+B\Pr)^{-1}}\\
    &\phantom{******}\times\vvert{((|v|^2-dh)^2-2h|v|^2)\Pr}\\
    &\leq C\nu^{-2}h^{-1}h^{2} = C\nu^{-2}h,
\end{split}
\ees
which proves
\be\label{mJp9}
((|v|^2-dh)^2-2h|v|^2)y\delta_y(h\alpha+B\Pr)^{-1}\Pr = \nu^{-2}O(h).
\ee
Combining \eqref{mJp8} and \eqref{mJp9} we have :
\be\label{eq:finalYY}
\nu^2Y^*Y(h\alpha+B\Pr)^{-1}\Pr=O(h).
\ee

Let us now study simultaneously the last two terms, using lemma \ref{lem:computprod} we have
\bes
\begin{split}
    &Y^*H_i(h\alpha+B\Pr)^{-1}\Pr = (-v_i(|v|^2-dh)\delta_{x_i}\delta_y+v_ihy\delta_{x_i})(h\alpha+B\Pr)^{-1}\Pr,\\
    &H_i^*Y(h\alpha+B\Pr)^{-1}\Pr = (-v_i(|v|^2-dh)\delta_{x_i}\delta_y+2v_ih\p_{x_i}V\delta_y)(h\alpha+B\Pr)^{-1}\Pr.
\end{split}
\ees
Thanks to Lemma \ref{lem:elliptreg} and since $\delta_x$ and $\delta_y$ commute
\bes
\begin{split}
    (\ast_3) &:= \vvert{\delta_{x_i}\delta_y(h\alpha+2d\nu^2hN_y)^{-1/2}(h\alpha+\Delta_{\frac{V}{2}})^{-1/2}}\\
    &\leq \vvert{\delta_y(h\alpha+2d\nu^2hN_y)^{-1/2}}\vvert{\delta_{x_i}(h\alpha+\Delta_{\frac{V}{2}})^{-1/2}}\\
    &\leq C\nu^{-1}h^{-\frac{1}{2}}.
\end{split}
\ees
Therefore, since $B$ commute with $\Delta_\frac{V}{2}$ and $N_y$, we get with Lemma \ref{lem:bornVPi} and Lemma \ref{lem:bornres}
\bes
\begin{split}
    (\ast_4) &:= \vvert{v_i(|v|^2-dh)\delta_{x_i}\delta_y(h\alpha+B\Pr)^{-1}\Pr}\\
    &\leq \vvert{\delta_{x_i}\delta_y(h\alpha+2d\nu^2hN_y)^{-1/2}(h\alpha+\Delta_{\frac{V}{2}})^{-1/2}}\vvert{v_i(|v|^2-dh)\Pr}\\
    &\phantom{*****}\times\vvert{(h\alpha+2d\nu^2hN_y)^{1/2}(h\alpha+B\Pr)^{-1/2}}\vvert{(h\alpha+\Delta_{\frac{V}{2}})^{1/2}(h\alpha+B\Pr)^{-1/2}}\\
    &\leq C\nu^{-1}h^{-\frac 12}h^{\frac 32}.
\end{split}
\ees
Consequently,
\bes
v_i(|v|^2-dh)\delta_y\delta_{x_i}(h\alpha+B\Pr)^{-1}\Pr = \nu^{-1}O(h).
\ees
Similarly, thanks to Lemma \ref{lem:elliptreg} and then Lemma \ref{lem:resolvestim}, considering $u\in L^2$
\bes
\begin{split}
    (\ast_5) &:= \vvert{y\delta_{x_i}(h\alpha+2d\nu^2hN_y)^{-1/2}(h\alpha+\Delta_{\frac{V}{2}})^{-1/2}u}^2\\
    &\leq \vvert{\delta_{x_i}(h\alpha+\Delta_{\frac{V}{2}})^{-1/2}}^2\vvert{y(h\alpha+2d\nu^2hN_y)^{-1/2}u}^2\\
    &\leq C\vvert{y(h\alpha+2d\nu^2hN_y)^{-1/2}u}^2\\
    &\leq C\< N_y(h\alpha+2d\nu^2hN_y)^{-1/2}u,(h\alpha+2d\nu^2hN_y)^{-1/2}u\>\\
    &\phantom{******}+Ch\vvert{(h\alpha+2d\nu^2hN_y)^{-1/2}u}^2\\
    &\leq C\big(\vvert{N_y(h\alpha+2d\nu^2hN_y)^{-1}}+h\vvert{(h\alpha+2d\nu^2hN_y)^{-1/2}}^2\big)\vvert{u}^2\\
    &\leq C(\nu^{-2}h^{-1}+hh^{-1}\alpha^{-1})\vvert{u}^2\\
    &\leq C\nu^{-2}h^{-1}(1+\nu^2h\alpha^{-1})\vvert{u}^2,
\end{split}
\ees
which leads to
\bes
\vvert{y\delta_{x_i}(h\alpha+2d\nu^2hN_y)^{-1/2}(h\alpha+\Delta_{\frac{V}{2}})^{-1/2}} \leq C\nu^{-1}h^{-\frac 12}(1+\nu\sqrt{h}\alpha^{-\frac12}).
\ees
And with very identical arguments, we get
\bes
\vvert{\p_{x_i}V\delta_y(h\alpha+2d\nu^2hN_y)^{-1/2}(h\alpha+\Delta_{\frac{V}{2}})^{-1/2}} \leq C\nu^{-1}h^{-\frac 12}(1+\alpha^{-\frac12}).
\ees
Therefore, with Lemma \ref{lem:bornVPi} and Lemma \ref{lem:bornres},
\bes
\vvert{v_ihy\delta_{x_i}(h\alpha+B\Pr)^{-1}}+\vvert{2v_ih\p_{x_i}V\delta_y(h\alpha+B\Pr)^{-1}}\leq C\nu^{-1}h(1+\alpha^{-\frac12}+\nu\sqrt{h}\alpha^{-\frac12}),
\ees
which leads to
\be\label{eq:finalYH_0H_0Y}
\nu(Y^*H_i+H_i^*Y)(h\alpha+B\Pr)^{-1}\Pr = (1+\alpha^{-\frac12}+\nu\sqrt{h}\alpha^{-\frac12})O(h).
\ee

Combining \eqref{eq:finalH_0H_0}, \eqref{eq:finalYY} and \eqref{eq:finalYH_0H_0Y} we have finally completely proved \eqref{eq:majJ1}.

To prove \eqref{eq:majJ2}, we first show that $\ooo A^*\Pr$ is bounded :
\bes
\begin{split}
    \ooo A^*\Pr &= \ooo Z\Pr{(h\alpha+B\Pr)^{-1}}\\
    \ooo Z\Pr u &= (\ooo (v\cdot\delta_x+\nu(|v|^2-dh)\delta_y)\Pr)u\\
    &=[\ooo,v\cdot\delta_x+\nu(|v|^2-dh)\delta_y]\Pr u\\
    &= -2dh^2\nu\delta_y\Pr u-2h^2(\delta_x+2\nu\delta_yv)\cdot\left(-\frac{v}{2h}\right)\Pr u\\
    &= h\left(v\cdot\delta_x+2\nu(|v|^2-dh)\delta_y\right)\Pr u.
\end{split}
\ees
Using Lemmas \ref{lem:bornVPi}, \ref{lem:bornres}, \ref{lem:resolvestim}, \ref{lem:elliptreg}, it proves
\bes
\ooo A^*\Pr = \alpha^{-\frac12}O(h).
\ees
Which leads to
\bes
\begin{split}
    \< A\ooo u,u\> &= \< \Pr A\ooo(1-\Pr) u,u\> = \<(1-\Pr)u,\ooo A^*\Pr u\>\\
    |\< A\ooo u,u\>| &\leq \vvert{\ooo A^*\Pr}\vvert{\Pr u}\vvert{(1-\Pr)u}\leq C\alpha^{-\frac12}h\vvert{\Pr u}\vvert{(1-\Pr)u}
\end{split}
\ees
which proves  \eqref{eq:majJ2}. Then, for \eqref{eq:majJ3} we show that $Z\Pr(h\alpha+B\Pr)^{-1/2}$ is bounded :
\bes
\begin{split}
    \vvert{Z\Pr(h\alpha+B\Pr)^{-1/2}u}^2 &= \< (Z\Pr)^*(Z\Pr)(h\alpha+B\Pr)^{-1/2}u,(h\alpha+B\Pr)^{-1/2}u\>\\
    &= h\vvert{u}^2-h^2\alpha\vvert{(h\alpha+B\Pr)^{-1/2}u}^2\\
\end{split}
\ees
which gives us
\bes
Z\Pr(h\alpha+B\Pr)^{-1/2} = O\big(\sqrt h\big).
\ees
Thanks to that result,
\bes
\begin{split}
    \< Zu, Au\> &= \< (1-\Pr)u,(\Pr Z)^*A(1-\Pr)u\>\\
    |\< Zu, Au\>| &\leq \vvert{(\Pr Z)^*A}\vvert{(1-\Pr)u}^2 = \vvert{-Z\Pr A}\vvert{(1-\Pr)u}^2\\
    &\leq \vvert{Z\Pr(h\alpha+B\Pr)^{-1}(Z\Pr)^*}\vvert{(1-\Pr)u}^2\\
    &\leq \vvert{Z\Pr(h\alpha+B\Pr)^{-1/2}}\vvert{(h\alpha+B\Pr)^{-1/2}(Z\Pr)^*}\vvert{(1-\Pr)u}^2\\
    &\leq \vvert{Z\Pr(h\alpha+B\Pr)^{-1/2}}^2\vvert{(1-\Pr)u}^2 \leq Ch\vvert{(1-\Pr)u}^2,
\end{split}
\ees
which completes the proof of Lemma \ref{lem:majorJprime}.

\ep

\subsection{Hypocercivity result.}
\begin{proposition}\label{prop:hypocoer}
For every $\delta_0>0$, there exists $C,h_0>0$ such that for all $h\in]0,h_0]$, $\gamma,\nu>0$ and for all $u\in D(P)\cap F_h^\bot$, one has
\bes
\Re\,\big\<Pu,(1+\delta(h)(A+A^*))u\big\>_{L^2}\geq C g(h) \Vert u\Vert_{L^2}^2
\ees
choosing $\alpha=\min(1,\nu^2h)$, where $\delta(h)=\delta_0\frac{g(h)}{h}$ and with $g(h)$ defined in \eqref{eq:g}.
\end{proposition}
\bp
For all $\delta>0$, let us define 
\bes
I_\delta=\Re\,\big\<P u,(1+\delta(A+A^*)u\big\>_{L^2}.
\ees
Using the decomposition $P=Z+\gamma\ooo$, and the skew-adjointness of $Z$ coming from \eqref{eq:adj0}, one gets
\bes
I_\delta=\gamma\<\ooo u,u\>+\delta\Re\<Pu,(A+A^*)u\>.
\ees
From the spectral properties of $\ooo$, it follows that
\be\label{eq:coer1}
I_\delta\geq \gamma h\Vert (1-\Pi_\rho)u\Vert^2+\delta\Re\<Pu,(A+A^*)u\>.
\ee
Denoting $J=\<Pu,(A+A^*)u\>$, one has
\bes
\begin{split}
J&=\<AZ u,u\>+\gamma\<A\ooo u, u\>+\<Z u,A u\>+\gamma\<\ooo u, Au\>,
\end{split}
\ees
and since $A=\Pi_\rho A$ and $\Pi_\rho \ooo=0$ it follows that 
$$
J=\<AZ \Pi_\rho u,u\>+J'
$$
with 
\bes
J'=\<AZ (1-\Pi_\rho )u,u\>+\gamma\<A\ooo u, u\>+\<Z u,A u\>.
\ees
Moreover, by definition, one has
$
AZ \Pi_\rho=(h\alpha+B\Pi_\rho)^{-1}hB\Pi_\rho
$.
Combined with Lemma \ref{eq:minorB} and taking $\alpha=\min(1,\nu^2h)$, this  shows that 
\be\label{eq:minorAZPi}
AZ \Pr \geq h\frac{c_0h\min(1,\nu^2h)}{h\min(1,\nu^2h)+c_0h\min(1,\nu^2h)} \Pr = \frac{c_0h}{1+c_0}\Pr = c_0'h\Pr.
\ee
Consequently, there exists $c_0'>0$ such that 
$$
\Re J\geq c_0'h\Vert \Pi_\rho u\Vert^2+\Re J'.
$$
Plugging this estimate into \eqref{eq:coer1} we get
\be\label{eq:coer2}
I_\delta\geq \gamma h\Vert (1-\Pi_\rho)u\Vert^2+\delta c_0'h\Vert \Pi_\rho u\Vert^2+\delta \Re J'.
\ee
Recall that 
\be\label{eq:defJprime}
J'=\<AZ(1-\Pi_\rho)u,u\>+\gamma\<A\ooo u, u\>+\<Z u,A u\>.
\ee
Combining \eqref{eq:coer2}, \eqref{eq:defJprime} and Lemma \ref{lem:majorJprime}, we get
\bes
\begin{split}
I_\delta &\geq -Ch \delta (1+(1+\nu\sqrt{h}+\gamma)\alpha^{-\frac12}) \Vert \Pi_\rho u\Vert\, \Vert(1- \Pi_\rho) u\Vert \\
&\phantom{******}+(\gamma h-Ch\delta)\Vert (1-\Pi_\rho)u\Vert^2+\delta c_0'h\Vert \Pi_\rho u\Vert^2\\
&\geq h\bigg(\gamma -C\delta-\frac{C^2\delta(1+(1+\nu\sqrt{h}+\gamma)\alpha^{-\frac12})^2}{2c_0'}\bigg)\Vert (1-\Pi_\rho)u\Vert^2+h\frac {\delta c_0'}2\Vert \Pi_\rho u\Vert^2.
\end{split}
\ees
Optimizing the right hand side  by taking 
\bes
\delta=\frac{2\gamma c_0'}{c_0'^2+2Cc_0'+C^2(1+(1+\nu\sqrt{h}+\gamma)\alpha^{-\frac12})^2},
\ees
we get
\bes
I_\delta\geq \frac {\delta hc_0'} 2\Vert u\Vert^2.
\ees

We shall say that $\delta(h) \asymp \tilde\delta(h)$ if there exists $C_1,C_2>0$ such that for $h$ small enough, $C_1\tilde\delta(h)\leq\delta(h)\leq C_2\tilde\delta(h)$. Therefore we have that

\bes
\delta \asymp \frac{\gamma}{1+(1+\nu^2h+\gamma^2)\alpha^{-1}}.
\ees

Recalling we took $\alpha=\min(1,\nu^2h)$, hence $\alpha^{-1}=\max(1,(\nu^2h)^{-1})$, one has

\underline{Case 1 :} $\nu^2h\leq1$

\bes
    \delta \asymp \frac{\gamma}{1+(1+\nu^2h+\gamma^2)(\nu^2h)^{-1}} \asymp \frac{\nu^2h\gamma}{1+\gamma^2} \asymp \min\big(\nu^2h\gamma,\frac{\nu^2h}\gamma\big).
\ees

\underline{Case 2 :} $1\leq\nu^2h$

\bes
    \delta \asymp \frac{\gamma}{1+\nu^2h+\gamma^2} \asymp \frac{\gamma}{\nu^2h+\gamma^2} \asymp \min\big(\frac{\gamma}{\nu^2h},\frac1\gamma\big).
\ees

This yields $\delta(h)\asymp\frac{g(h)}{h}$ and therefore

\bes
I_\delta\geq Cg(h)\Vert u\Vert^2
\ees
for some new constant $C>0$ independent of $h,\gamma$ and $\nu$. This proves the proposition.

\ep

\subsection{Proof of Theorem \ref{thm}.}
Let $u\in D(P)\cap F_h^\bot$. On the one hand, with a mere Cauchy-Schwartz, we have
\be\label{eq:thm1}
\Re\< (P-z)u, (1+\delta(h)(A+A^*))u \> \leq \vvert{(P-z)u}\vvert{1+\delta(h)(A+A^*)}\vvert{u}.
\ee
On the other hand, we can see thanks to Proposition \ref{prop:hypocoer} that
\be\label{eq:eq00}
\Re\< (P-z)u, (1+\delta(h)(A+A^*))u \> \geq Cg(h)\vvert{u}^2-\Re(z\< u,(1+\delta(h)(A+A^*))u\>).
\ee
Using that for $\delta_0$ small enough,
\bes
\vvert{\delta(h)(A+A^*)} \leq 2\delta(h)\alpha^{-\frac12} < 1,
\ees
we have that $1+\delta(h)(A+A^*)$ is positive and $\vvert{1 + \delta(h)(A+A^*)}\leq2$, hence we can simplify \eqref{eq:eq00} to
\be\label{eq:thm2p}
\Re\< (P-z)u, (1+\delta(h)(A+A^*)u \> \geq Cg(h)\vvert{u}^2-2\, |\Re(z)|\vvert{u}^2.
\ee
Moreover \eqref{eq:thm1} becomes
\be\label{eq:thm1p}
\Re\< (P-z)u, (1+\delta(h)(A+A^*))u \> \leq 2\vvert{(P-z)u}\vvert{u}.
\ee
We can combine \eqref{eq:thm1p} and \eqref{eq:thm2p} to have
\bes
\vvert{(P-z)u} \geq \frac C2g(h)\vvert{u}-|\Re(z)|\vvert{u}.
\ees
Finally, taking $0<c_0<\frac C2$ and noting $c_1=\frac C2-c_0>0$, we have for $|\Re z|\leq c_0g(h)$
\be\label{eq:PmzuFbot}
\vvert{(P-z)u} \geq c_1g(h)\vvert{u}.
\ee

And we can now deduce the second part of Theorem \ref{thm} from that, following the same sketch of proof as in \cite{No23}.

Let $\m\in\uuu^{(0)}$, by recalling $f_\m(x,v,y)=\chi_\m(x) e^{-(f(x,v,y)-f(\m))/h}$, as we know from \eqref{eq:H_0YNmu} that $e^{-f/h}\in\Ker\ooo\cap\Ker Y$, we obtain :
\bes
P(f_\m)=H_0(f_\m)=hv\cdot\nabla\chi_\m e^{-(f-f(\m))/h}=O(h^{-\frac d2 + \frac34}e^{-c_\m/h})
\ees
with $c_\m=\inf_{\supp\nabla\chi_\m} f-f(\m)>0$ (because $\chi_\m$ is constant of order $h^{-\frac d2 - \frac14}$ near $\m$).

Moreover, since the $(f_\m)_{\m\in\uuu^{(0)}}$ are orthonormal, we actually have
\be\label{eq:Pfm}
\forall u\in F_h,\ \vvert{Pu}=O(h^{-\frac d2 + \frac34}e^{-c_f/h})\vvert{u}
\ee
where $c_f=\min_{\m\in\uuu^{(0)}}c_\m>0$. Furthermore, \eqref{eq:Pfm} is also true replacing $P$ by $P^*$ because $P^*(f_\m)=-H_0(f_\m)$. And because
\bes
P^*Pf_\m = -H_0(hv\cdot\chi_\m)e^{-(f-f(\m))/h} = h^2(\nabla V\cdot\nabla\chi_\m - \Hess\chi_\m v\cdot v)e^{-(f-f(\m))/h},
\ees
\eqref{eq:Pfm} is still valid replacing $P$ by $P^*P$.

We denote by $\Pi$ the orthogonal projector on $F_h$. Let $z\in\{|\Re z|\leq c_0g(h)\}$ such that $|z|\geq c_0'g(h)$ with $0<c_0'\leq c_1$, and $u\in D(P)$
\bes
\begin{split}
    \vvert{(P-z)u}^2 &= \vvert{(P-z)(\Pi+\Id-\Pi)u}^2\\
    &= \vvert{(P-z)(\Id-\Pi)u}^2+\vvert{(P-z)\Pi u}^2\\&\hspace{1.5cm}+2\Re\<(P-z)(\Id-\Pi)u,(P-z)\Pi u\>.
\end{split}
\ees
One has
\bes
\begin{split}
    &\vvert{(P-z)(\Id-\Pi)u}^2 \geq c_1^2g(h)^2\vvert{(\Id-\Pi)u}^2,\\
    &\vvert{(P-z)\Pi u}^2 \geq (\vvert{P\Pi u}-\vvert{z\Pi u})^2 \geq \vvert{z\Pi u}(\vvert{z\Pi u}-2\vvert{P\Pi u}).
\end{split}
\ees
Using that $|z|\geq c_0'g(h)\geq c_0'e^{-c/(2h)}\geq c_0'e^{-c/h}$ with $c<c_f$ thanks to \eqref{eq:gthm2}, we have using \eqref{eq:Pfm}
\bes
\vvert{(P-z)\Pi u}^2 \geq \frac{|z|^2}{2}\vvert{\Pi u}^2.
\ees
We can also see, studying each term in the scalar product :
\bes
\begin{split}
    (\ast_6) :\negthinspace&=\Re\<(P-z)(\Id-\Pi)u,(P-z)\Pi u\>\\
    &= \Re\big(\<P(\Id-\Pi)u,P\Pi u\>-z\<(\Id-\Pi)u,P\Pi u\> - \bar{z}\<P(\Id-\Pi)u,\Pi u\>\big)\\
    &\leq (1+|z|)\vvert{(\Id-\Pi)u}\vvert{\Pi u}O(h^{-\frac d2 + \frac34}e^{-c_f/h})\\
    &\leq \Big(\vvert{u}^2+|z|^2\vvert{\Pi u}^2+\vvert{(\Id-\Pi)u}^2\Big)O(h^{-\frac d2 + \frac34}e^{-c_f/h}),
\end{split}
\ees
hence
\bes
\begin{split}
    \vvert{(P-z)u}^2 &\geq c_1^2g(h)^2\vvert{(\Id-\Pi)u}^2+\frac{|z|^2}3\vvert{\Pi u}^2\\
    &\phantom{*******}+(\vvert{u}^2+\vvert{(\Id-\Pi)u}^2)O(h^{-\frac d2 + \frac34}e^{-c_f/h})\\
    &\geq \frac{c_0'^2}3g(h)^2\vvert{u}^2+(\vvert{u}^2+\vvert{(\Id-\Pi)u}^2)O(h^{-\frac d2 + \frac34}e^{-c_f/h})\\
    &\geq \frac{c_0'^2}{4}g(h)^2\vvert{u}^2
\end{split}
\ees
for $h$ small enough, using \eqref{eq:gthm2}. It leads to
\be\label{eq:injectivePz}
\vvert{(P-z)u} \geq \frac{c_0'}{2}g(h)\vvert{u}.
\ee

By using the same arguments for $P^*$ we have the same result for it (the key point is that $e^{-f/h}$ is in the kernel of $\ooo,H_0$ and $Y$ hence it also is in $P^*$'s one). It just remains to show that $P-z$ is surjective in order to obtain the resolvent estimate, we show it the usual way, by showing that $\Ran(P-z)$ is closed and dense.

Let $u_n\in D(P)$ and $v\in L^2$ such that $(P-z)u_n\to v$, therefore $((P-z)u_n)_{n\in\N}$ is Cauchy and so is $(u_n)_{n\in\N}$ thanks to \eqref{eq:injectivePz}, hence there exists $u\in L^2$ such that $u_n\to u$. Because the convergence is also true in $\ddd'$, we have that $(P-z)u=v$ in $\ddd'$, and since $v\in L^2$, so is $(P-z)u$, thus $u\in D(P)$ and $\Ran(P-z)$ is closed. Now to show that $\Ran(P-z)$ is dense, we use \eqref{eq:injectivePz} for $P^*$ and so $\Ker(P^*-\overline{z})=\{0\}$.

All this leads to the resolvent estimate :
\be\label{resol_esti}
\vvert{(P-z)^{-1}}\leq\frac{2}{c_0'g(h)}.
\ee

Hence, $P$ has no spectrum in $\{|\Re z|\leq c_0g(h)\}\cap\{|z|\geq c_0'g(h)\}$. Thanks to Proposition \ref{prop:accretive}, we know that $P$ is maximally accretive and therefore $P - z$ is invertible for all $\Re z < 0$. Moreover we easily see that for $u\in D(P)$,
\bes
\vvert{(P-z)u}\vvert{u} \geq \Re \<(P-z)u,u\> \geq -\Re z\vvert{u}^2,
\ees
thus for all $\Re z < 0$,
\bes
\vvert{(P-z)^{-1}} \leq \frac{1}{-\Re z},
\ees
which extends \eqref{resol_esti} :
\bes
\forall z\in\{\Re z\leq c_0g(h)\}\cap\{|z|\geq c_0'g(h)\},\ \ \vvert{(P-z)^{-1}}\leq\frac{2}{c_0'g(h)}.
\ees

We will show that on $\{\Re z\leq c_0g(h)\}$, $P$ has exactly $n_0=\dim F_h$ eigenvalues counted with multiplicities. By denoting $D =D(0,c_0' g(h))$ the disk in $\C$ centered at $0$ of radius $c_0' g(h)$, let us denote
\bes
\Pi_0=\frac{1}{2i\pi}\int_{\p\negthinspace D}(z-P)^{-1}dz
\ees
the Riesz projector on the small eigenvalues. We start by proving the following lemma
\begin{lemma}\label{lem:Pi_0}
We have $\vvert{P\Pi_0}\leq 2c_0' g(h)$.
\end{lemma}

\bp
\bes
P\Pi_0 = \frac{1}{2i\pi}\int_{\p\negthinspace D}P(z-P)^{-1}dz = \frac{1}{2i\pi}\int_{\p\negthinspace D}z(z-P)^{-1}dz,
\ees
hence
\bes
\vvert{P\Pi_0}\leq (c_0' g(h))^2\frac{2}{c_0'g(h)} = 2c_0' g(h)
\ees
thanks to \eqref{resol_esti}.

\ep

We first prove that $\dim\Ran\Pi_0\leq n_0$.

By contradiction, let us suppose $F_h^\bot\cap\Ran\Pi_0\neq\emptyset$ and so let us take $u\in F_h^\bot\cap\Ran\Pi_0$ of norm one. Since $u\in\Ran\Pi_0$, by Lemma \ref{lem:Pi_0}, $\vvert{Pu}\leq 2c_0' g(h)$, but because $u\in F_h^\bot$ we can use \eqref{eq:PmzuFbot} and so $\vvert{Pu}\geq c_1g(h)$. Taking $c_0'$ low enough, we have the contradiction we aimed for and thus, $\dim\Ran\Pi_0\leq n_0$.

For the converse inequality, we have
\bes
\Pi_0-\Id=\frac{1}{2i\pi}\int_{\p\negthinspace D}z^{-1}(z-P)^{-1}Pdz
\ees
and therefore
\be\label{quasiortho}
\begin{split}
    \eps_\m&=\Pi_0f_\m-f_\m=\frac{1}{2i\pi}\int_{\p\negthinspace D}\underbrace{(z-P)^{-1}}_{= O(g(h)^{-1})}\underbrace{P(f_\m)}_{= O(h^{-\frac d2 + \frac34}e^{-c_f/h})}\frac{dz}z\\
    &= O(g(h)^{-1}h^{-\frac d2 + \frac34}e^{-c_f/h}) = O(e^{-\frac{c_f}{2h}})
\end{split}
\ee
using \eqref{resol_esti} and \eqref{eq:Pfm} and the hypothesis \eqref{eq:gthm2}.

Let us suppose $\sum_{\m\in\uuu^{(0)}} a_\m\Pi_0 f_\m=0$ with $\sum_{\m\in\uuu^{(0)}}  |a_\m|^2=1$, since $\Pi_0 f_\m=f_\m+\eps_\m$, we have for all $\m'\in\uuu^{(0)}\sum_{\m\in\uuu^{(0)}}  a_\m(\delta_{\m,\m'}+\<\eps_\m,f_{\m'}\>)=0$ and so for all $\m$, $a_\m=O(e^{-c/h})$ for some $c<\frac{c_f}{2}$, which is in contradiction with $\sum |a_\m|^2=1$. We deduce that $\dim\Ran\Pi_0\geq n_0$ and hence, with what we already showed, $\dim\Ran\Pi_0= n_0$.

And so we can say we have
\bes
\sigma(P)\cap\{\Re z\leq c_0g(h)\}=\{\lambda_\m(h),\m\in\uuu^{(0)}\}\subset D(0,\frac{c_1}{2}g(h)).
\ees

It only remains to show that $\lambda_\m(h)=O(e^{-c/h})$ for some rather explicit $c>0$. Using the definition of $\Pi_0$ and \eqref{resol_esti}, we have that $\vvert{\Pi_0}\leq2$, moreover noticing that $\Ran\Pi_0$ is $P$-stable and that $(\Pi_0f_\m)_{\m\in\uuu^{(0)}}$ is one of its basis,
\bes
\vvert{P\Pi_0 f_\m} = \vvert{\Pi_0 Pf_\m} \leq 2\vvert{Pf_\m} = O(h^{-\frac d2 + \frac34}e^{-c_f/h}).
\ees
Therefore, $P_{|\Ran\Pi_0}=O(h^{-\frac d2 + \frac34}e^{-c_f/h})$ hence $\sigma(P_{|\Ran\Pi_0})\subset D(0,Ch^{-\frac d2 + \frac34}e^{-c_f/h})$ for some $C>0$. In other words, we have $|\lambda_\m(h)|\leq e^{-(c_f-\eps)/h}$ for all $\eps>0$.

\section{Sharp quasimodes}\label{sec:LocalQuasimode}

We now want to have a better view on the small eigenvalues of $P$. For this purpose, we are going to build sharp quasimodes, and so we are following the steps of \cite[section 3\&4]{BoLePMi22}. Their theorem does not apply here because our operator $P$ does not satisfy the hypothesis they labeled (Harmo) as explained in the last paragraph before the statements. We therefore have to rewrite the proof using tricks to avoid that necessity. As in this reference, given $\s\in\uuu^{(1)}$ and $U$ a neighborhood of $\s$ we look for an approximate solution to the equation $P\tilde u=0$, with $\tilde u\in\ccc^\infty (\overline U)$ under the form
$$\tilde u=ue^{-(f-f(\m))/h}$$
where we recall $e^{-(f-f(\m))/h}\in\Ker P$ with $f(x,v,y)=\frac{V(x)}{2}+\frac{|v|^2+y^2}{4}$. And we set
\bes
u(x,v,y)=\int_0^{\ell(x,v,y,h)}\zeta(s/\tau)e^{-s^2/2h}ds
\ees
where the function $\ell\in\ccc^\infty(U)$ has a classical expansion $\ell\sim\sum_{j\geq0}h^j\ell_j$ in $\ccc^\infty(U)$. Here, $\zeta$ denotes a fixed smooth even function equal to $1$ on $[-1,1]$ and supported in $[-2,2]$, and $\tau > 0$ is a small parameter which will be fixed later. The object of this section is to construct the function $\ell$. In the following, we will use $X$ instead of $(x,v,y)$ to simplify the equations.

We see that our operator $P$ can be written as in \cite{BoLePMi22} : 
\be\label{eq:formP}
P=-h\div\circ A\circ h\nabla +\frac{1}{2}(b\cdot h\nabla + h\div\circ b)+c
\ee
with
\bes
A=\begin{pmatrix}0&0&0\\0&\gamma\Id&0\\0&0&0\end{pmatrix},\
b=\begin{pmatrix}v\\-\p_xV-\nu yv\\\nu(|v|^2-dh)\end{pmatrix} \mbox{ and }
c=\gamma\Big(\frac{|v|^2}{4}-h\frac d2\Big).
\ees

Since $P$ is of the form \eqref{eq:formP}, we can apply \cite[Lemma 3.1]{BoLePMi22}, and we get
\be\label{eq:Pwr}
P(ue^{-(f-f(\m))/h})=h(w+r)e^{-\big(f-f(\m)+\frac{\ell^2}{2}\big)/h}
\ee
with $r$ vanishing around $\s$ and
\be\label{eq:w}
w = (b+2A\nabla f)\cdot\nabla\ell+\ell A\nabla\ell\cdot\nabla\ell-h\div A\nabla\ell.
\ee
In the following, we will consider $\gamma,\nu>0$ fixed and $\s=0$. Under these hypotheses, $w$ can be expressed in powers of $h$, $w\sim \sum_{j\geq0}h^jw_j$.

Foreshadowing the suitable estimates we will need in the end, we want to solve $w=O(X^4+hX^2+h^2)$ in order to have
\bes
\vvert{P(ue^{-(f-f(\m))/h})} = O(h^2)\sqrt{\lambda_\m},
\ees
this order is the lowest that will give us precise results on the low lying eigenvalues $\lambda_\m$, hence this choice. Thus we decide to take $\ell=\ell_0+h\ell_1$, which gives us $w=w_0+hw_1+O(h^2)$ with
\be\label{eq:w0w1}
\left\{\begin{aligned}
    &w_0 = (b^0+2A\nabla f)\cdot\nabla\ell_0+A\nabla\ell_0\cdot\nabla\ell_0\ \ell_0,\\
    &w_1 = (b^0+2A(\nabla f+\ell_0\nabla\ell_0))\cdot\nabla\ell_1+A\nabla\ell_0\cdot\nabla\ell_0\ \ell_1+R_1,
\end{aligned}\right.
\ee
where $b= b^0+hb^1$ and $R_1 = b^1\cdot\nabla\ell_0-\div A\nabla\ell_0$. As in \cite{BoLePMi22} we call eikonal equation $w_0=0$ and transport equation $w_1=0$, we now are going to solve the eikonal equation up to the fourth order, and the transport equation up to the second one.

\subsection{Solving the eikonal equation.}
Unlike in \cite{BoLePMi22}, the outgoing manifolds of the flow passing through the saddle point are not a Lagrangian ones that project nicely on the $X$-space. Therefore, we need to find an other way to solve that equation, we will consider homogeneous polynomials to simplify it following \cite[Remark 2.3.9]{He88_01}.

We introduce $\ppp^j_{hom}$ the set of homogeneous polynomial of degree $j$ and we consider
\be\label{eq:ellpol}\ell_0=\sum_{j=0}^3\ell_{0,j},
\ee
with $\ell_{0,j}\in\ppp^j_{hom}$. In the following, we will need to have $\ell_{0}(\s) = \ell_0(0) = 0$ (recalling we assumed $\s=0$), therefore we need to set $\ell_{0,0}=0$. We also denote $\ell_{0,1}(X)=\xi\cdot X$ for a certain $\xi\in\R^{2d+1}$ to be determined. Thanks to \eqref{eq:ellpol} and \eqref{eq:w0w1}, we have that $w_0$ also has a similar development $w_0=\sum_{j=0}^3w_{0,j}+O(X^4)$, $w_{0,j}\in\ppp^j_{hom}$ with $w_{0,0} = A\xi\cdot\xi \ \ell_{0,0} = 0$.

As in \cite{BoLePMi22}, we denote by $H$ and $B$ the matrices $\Hess_\s f$ and $db^0(\s)$ respectively. We also denote
$$\Lambda=2HA+B^T=\begin{pmatrix}
0&-\Hess_\s V&0\\
\Id&\gamma\Id&0\\
0&0&0
\end{pmatrix}.
$$
We first want to study the spectrum of $\Lambda$. To that purpose, let us consider $\Lambda' = \begin{pmatrix}0&-\Hess_\s V\\\Id&\gamma\Id\end{pmatrix}$ which is obtained when considering the usual Kramers-Fokker-Planck equation derived from the Langevin dynamics
\bes
\left\{\begin{aligned}
    dx_t &= v_tdt,\\
    dv_t &= -\nabla V(x_t)dt -\gamma v_tdt +\sqrt{2\gamma h}dB_t.
\end{aligned}\right.
\ees
This equation and its generator are non-degenerated, thus we can apply \cite{BoLePMi22}, in particular we have \cite[Lemma 1.4]{BoLePMi22} which states that $\Lambda'$ has one negative eigenvalue we call $-\mu$ and the rest of its spectrum is included in $\{\Re z >0\}$.

Now noticing that at the order $1$, we have $\nabla f(X)\sim HX$ and $b^0(X)\sim BX$, we have that $w_{0,1}=0$ becomes
\be\label{eik2}
\Lambda\xi\cdot X+(A\xi\cdot\xi)\xi\cdot X=0,
\ee
with unknown $\xi$ and must be true for any $X$ near $0$. Taking $\xi$ an eigenvector of $\Lambda$ associated with the negative eigenvalue $-\mu$ will solve the equation, we will just have to chose the right vector on $\R\xi$.

A good choice is to consider $\xi = t\begin{pmatrix}\alpha e_1\\e_1\\0\end{pmatrix}$ with $t\in\R$, $e_1$ and eigenvector of $\Hess_\s V$ associated to its only negative eigenvalue $-\eta<0$ (because $\s$ is a critical point of index $1$) and $\alpha$ to be determined. Therefore, solving $\Lambda\xi = -\mu\xi$ leads to $\alpha = -\frac12(\gamma + \sqrt{\gamma^2 + 4\eta})$ and 
\be\label{eq:mu}
\mu=\frac{1}{2}(-\gamma+\sqrt{\gamma^2+4\eta})>0,
\ee

Going back to \eqref{eik2}, $t$ is determined so that
\be\label{eq:Axiximu}
(A\xi\cdot\xi)-\mu=0,
\ee
we thus have $t = \frac{1}{|e_1|}\sqrt{\frac{\mu}{\gamma}}$ (with $e_1\neq0$ because it is an eigenvector).

We now notice that for $j\in\{2,3\}$, we have
\be\label{eq:w0j}
w_{0,j}=\lll\ell_{0,j}+R_{0,j},
\ee
with $\lll = \Upsilon X\cdot\nabla +\mu$ an endomorphism of $\ppp^j_{hom}$, $\Upsilon = \Lambda^T+2A\Pi_\xi$, $\Pi_\xi(X) = X\cdot\xi \ \xi$ and where $R_{0,j}$ is a smooth function of $\ell_{0,k}$ and $\nabla\ell_{0,k}$ for $k<j$.

Solving \eqref{eq:w0j} by homogeneous polynomial is a technique we take from \cite[Chapter 3]{DiSj99_01}, although we will not solve it up to $O(X^\infty)$.

In a basis of $\R^{2d+1}$ adapted to $\xi$ in which $\Lambda$ is upper triangular, we have that only the first entry of its diagonal has negative real part, being $-\mu$. Moreover in that same basis, $2\Pi_\xi A$ has zeros outside its first row, and the first element of that row is $2\mu$. Hence all the eigenvalues of $\Upsilon^T$, and thus of $\Upsilon$, have non-negative real part.

In other words, we have $\sigma(\Upsilon)\subset \{\Re z\geq0\}$, then $\sigma(\Upsilon X\cdot\nabla)\subset \{\Re z\geq0\}$ using \cite[Lemma A.1]{BoLePMi22} (in this lemma, the authors only consider $\{\Re z>0\}$, but there is no difficulty expanding the result to $\{\Re z\geq0\}$ either by continuity or just by doing the same proof). Because $\mu>0$, we thus have that $\lll$ is invertible and so we can solve $w_0=O(X^4)$.

\subsection{Solving the transport equation.}
The transport equation is much simpler to solve after having solved the eikonal one. Taking $\ell_1 = \ell_{1,0} + \ell_{1,1}$ with $\ell_{1,j}\in\ppp^j_{hom}$, we have $w_1=w_{1,0}+w_{1,1}+O(X^2)$, $w_{1,j}\in\ppp^j_{hom}$ and
\bes\left\{
\begin{aligned}
&w_{1,0} = \mu\ell_{1,0}+b^1\cdot\xi-\div A\nabla\ell_{0,2},\\
&w_{1,1} = \lll\ell_{1,1}+R_{1,1},
\end{aligned}\right.
\ees
with $R_{1,1}$ a smooth function of $\ell_{0}$, $\ell_{1,0}$ and their derivatives up to the second order. The first equation is easily solved ($\mu\neq0$) and the second is solved using the same method as for $w_{0,2}$ and $w_{0,3}$.

After this we now have 
\be\label{eq:Pchi}
P(ue^{-(f-f(\m))/h})=hO(X^4+hX^2+h^2)e^{-\big(f-f(\m)+\frac{\ell^2}{2}\big)/h}.
\ee

In the next section, we will perform some Laplace methods applied to the right hand side of \eqref{eq:Pchi}, hence we need to determine that $\s$ is a local minimum of $f+\frac{\ell_0^2}{2}$. We thus have to state \cite[Lemma 3.3]{BoLePMi22}'s result and adapt \cite[Lemma 4.1]{BoLePMi22}'s proof as we don't have their Lagrangian manifold and its generating function which directly proves this lemma, let us recall the result :

\begin{lemma}\label{lem:xiortho}
We have
\be\label{eq:detHplusPixi}
\det\Hess_\s\big(f+\frac{\ell_0^2}{2}\big) = -\det H.
\ee
and hence, recalling that $\s\in\uuu^{(1)}$
\bes
\Hess_\s\big(f+\frac{\ell^2_0}2\big)>0.
\ees
Thus, around $\s$,
\be\label{eq:skindamin}
X-\s\in\xi^\bot\implies f(X)>f(\s).
\ee
\end{lemma}

\bp
We first observe that  :
\bes
\Hess_\s\big(f+\frac{\ell_0^2}{2}\big)=H+\Pi_\xi.
\ees
We thus have that \eqref{eq:detHplusPixi} is equivalent to
\bes
\det E=-1,
\ees
where $E=\Id+H^{-1}\Pi_\xi$. We first observe that $\xi^\bot$ is stable by $E$ and that $E_{|\xi^\bot}=\Id$. On the other hand, one has
\bes
E\xi\cdot\xi = \vvert{\xi}^2(1+H^{-1}\xi\cdot\xi).
\ees
But, $H(2A+H^{-1}B^T)\xi=\Lambda\xi=-\mu\xi$ gives $(2A+H^{-1}B^T)\xi\cdot\xi=-\mu H^{-1}\xi\cdot\xi$.

Looking at the skew-adjoint part of $P$ and using \eqref{eq:PPstarexp}, we obtain
\bes
h\div b - 2b\cdot\nabla f = 0,
\ees
identifying the first term in the classical expansion we get
\bes
b^0\cdot\nabla f = 0,
\ees
and knowing that $\nabla f(X) \sim HX$ and $b^0(X)\sim BX$ it follows that $B^TH$ is antisymmetric and so is $H^{-1}B^T=H^{-1}(B^TH)H^{-1}$ because $H$ is symmetric. Hence we have
\bes
H^{-1}\xi\cdot\xi = -\frac2\mu A\xi\cdot\xi = -2,
\ees
using \eqref{eq:Axiximu}, which leads to $E\xi\cdot\xi = -\vvert{\xi}^2$. Taking a basis adapted to $\xi^\bot$ completed with $\xi$ we can easily compute $\det E$ and obtain the aimed result.

For \eqref{eq:skindamin}, around $\s$,
\bes
f(X)=f(\s)+\frac{1}{2}H(X-\s)\cdot(X-\s)+O(|X-\s|^3).
\ees
But, over $\xi^\bot$, $H=H+\Pi_\xi=\Hess_\s\big(f+\frac{\ell_0^2}{2}\big)>0$ by what we have done just before.

\ep

\begin{remark}\label{rem:signell}
    Notice that the function $-\ell$ solves the equations the same way $\ell$ does. For now, the sign does not matter, but we will fix it in the next section when properly constructing the quasimodes.
\end{remark}

\section{Simpler case of a two wells function}\label{sec:TwoWells}

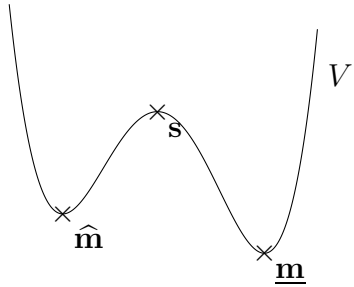
\begin{figure}[h]
    \centering
    \begin{tikzpicture}[scale=1]
    \draw (-1.1874,-1.3517) node{$\times$} node[below right]{$\widehat{\m}$};
    \draw (0.0706,.0089) node{$\times$} node[below right]{$\s$};
    \draw (1.4918,-1.8753) node{$\times$} node[below right]{$\underline{\m}$};
    \draw (2.5,.5) node{$V$};
    \draw plot [domain=-0.95:1.1,samples=1000] (\x*2,8*\x*\x*\x*\x-2*\x*\x*\x-7*\x*\x+\x/2);
    \end{tikzpicture}
    \caption{Representation of a typical two-wells Morse function}
    \label{fig:two_wells}
\end{figure}

In the following, we will restrain our study to a potential $V$ being a two-wells function, in other words, a function satisfying $\uuu^{(0)}=\{\underline{\m},\widehat{\m}\}$ and $\uuu^{(1)}=\{\s\}$ (see fig. \ref{fig:two_wells}). Moreover we assume that the wells have different depths, and we choose $\underline{\m}$ to be the deepest, namely $V(\underline{\m}) < V(\widehat{\m})$. It is a much simpler case that allows us not to consider many geometric constructions but still have interesting results and see how the $\ell$ we built is useful in these findings, see \cite[Definition 1.3]{BoLePMi22} for a complete description of the geometric construction on a more general case.

In this configuration we know from Theorem \ref{thm} that $P$ has exactly two low lying eigenvalues, among which $0$ that we decide to associate to $\underline{\m}$, the other one still not precisely known is associated to $\widehat{\m}$. This choice will appear to be relevant later on, when the exact form of the eigenvalue will be explicit. This is why we consider two wells and so we will focus on the second smallest eigenvalue, the aim of this section is to prove Theorem \ref{thm3}.

We define several sets following the description of \cite[Section 4]{BoLePMi22}, for $\tau,\delta>0$ :
\bes
\begin{split}
   &\bbb_{\tau,\delta}=\{f\leq f(\s)+\delta\}\cap\{X\in\R^{2d+1},\ |\xi\cdot(X-\s)|\leq\tau\},\\
   &E_{\tau,\delta}=\{f\leq f(\s)+\delta\}\setminus\ccc_{\tau,\delta},
\end{split}
\ees
where $\ccc_{\tau,\delta}$ denotes the connected component of $\bbb_{\tau,\delta}$ containing $\s$. We note $E^+_{\tau,\delta}$ the connected component of $E_{\tau,\delta}$ containing $\widehat{\m}$ and $E^-_{\tau,\delta}$ its complement in $E_{\tau,\delta}$. One can show that for $\tau_0,\delta_0$ small enough, for every $\tau\in\,]0,\tau_0],\delta\in\,]0,\delta_0]$, we have $\underline{\m}\in E^-_{\tau,\delta}$.

These are useful to define the following cutoff properly
\bes
\chi_\ell(X)=\left\{\begin{array}{cc}
     +1& \mbox{for }X\in E^+_{4\tau,4\delta}, \\
     -1& \mbox{for }X\in E^-_{4\tau,4\delta}
\end{array}\right.
\ees
and
\bes
\chi_\ell(X)=C_h^{-1}\int_0^{\ell(X)}\zeta(r/\tau)e^{-\frac{r^2}{2h}}dr\quad \mbox{for }X\in\ccc_{4\tau,4\delta},
\ees

where $C_h=\frac 12\int_{-\infty}^{+\infty}\zeta(r/\tau)e^{-\frac{r^2}{2h}}dr$, $\zeta\in\ccc^\infty(\R,[0,1])$ is even and satisfies $\zeta=1$ on $[-1,1]$ and $\zeta=0$ outside $[-2,2]$. And here, $\ell$ is the function built in the previous section which sign (see remark \ref{rem:signell}) is chosen so that there exists a neighborhood $\www$ of $\s$ such that $E^+_{4\tau,4\delta}\cap\www$ is included in $\{\xi\cdot(X-\s)>0\}$.

We notice by the way that $\exists\beta>0, C_h^{-1}=\sqrt{\frac{2}{\pi h}}(1+O(e^{-\beta/h}))$. For $h,\delta>0$ small, $\ell\sim\ell_0\sim \xi\cdot(X-\s)=\pm4\tau$ on the boundary of $\ccc_{4\tau,4\delta}$ within $\{f\leq f(\s) + 4\delta\}$. But $\zeta(\cdot/\tau)$ is vanishing outside $(-2\tau,2\tau)$, therefore $\chi_\ell$ is indeed a smooth function on $\{f\leq f(\s) + 4\delta\}$. To have a cutoff defined properly on $\R^d$, we introduce
\bes
\theta(X)=\left\{\begin{array}{cl}
     1& \mbox{for }X\in \{f\leq f(\s)+\delta\}, \\
     0& \mbox{for }X\in \R^d\setminus\{f\leq f(\s)+2\delta\}
\end{array}\right.
\ees
and smooth between these sets.

Hence, we have $\theta \chi_\ell\in\ccc_c^\infty(\R^{2d+1},[-1,1])$ and $\supp \theta \chi_\ell\subset\{f\leq f(\s)+2\delta\}$.

\begin{defin}For $\tau>0$ and then $\delta,h>0$ small enough, we define the quasimodes 
\bes
\left\{\begin{array}{l}
    \psi_{\underline{\m}}(X)=2e^{-\frac{f(X)-f(\underline{\m})}{h}},\\
    \psi_{\widehat{\m}}(X)=\theta(X)(\chi_\ell(X)+1)e^{-\frac{f(X)-f(\widehat{\m})}{h}}.
\end{array}\right.
\ees
And in the same time, we define the normalized quasimodes
\bes
\phii_{\underline{\m}}=\frac{\psi_{\underline{\m}}}{\vvert{\psi_{\underline{\m}}}},\quad \quad\phii_{\widehat{\m}}=\frac{\psi_{\widehat{\m}}}{\vvert{\psi_{\widehat{\m}}}}.
\ees
\end{defin}



For shortness, we write $D_{X^*}=|\det\Hess_{X^*}(f)|^{1/2}$ for $X^*\in\uuu$. We recall that $\gamma,\nu>0$ are fixed.

\begin{proposition}\label{prop:interactionmatrix}
For $\tau>0$ and then $\delta>0$ small enough, there exists $C>0$ such that for every $\m,\m'\in\uuu^{(0)}=\{\widehat{\m},\underline{\m}\}$ and $h>0$ small,
\begin{itemize}
    \item[$i)$] $\<\phii_\m,\phii_{\m'}\>=\delta_{\m,\m'}+O(e^{-C/h}),$
    \item[$ii)$] $\<P\phii_\m,\phii_\m\>=he^{-2\tilde S(\m)/h}\frac{\mu(\s)}{2\pi}\frac{D_\m}{D_\s}(1+O(h))$,
    \item[$iii)$] $\vvert{P\phii_\m}^2=O(h^4)\<P\phii_\m,\phii_\m\>$,
    \item[$iv)$] $\vvert{P^*\phii_\m}^2=O(h)\<P\phii_\m,\phii_\m\>$,
\end{itemize}
where $\tilde S(\widehat{\m})=f(\s)-f(\widehat{\m})$ and $\tilde S(\underline{\m})=+\infty$, and with $\mu(\s)$ defined in \eqref{eq:mu}.
\end{proposition}
\begin{remark}
We have built our constant $\mu$ to be positive and thus $-\mu$ is the negative eigenvalue of $\Lambda$ while \cite{BoLePMi22} did it the other way, hence the lack of absolute value in $ii)$.
\end{remark}

\bp
Noticing $f$ attains its absolute minimum at $\m$ on $\supp\psi_\m$, by using a Laplace method on $\psi_\m$ we obtain
\be\label{eclpsi}
\vvert{\psi_\m}=2(\pi h)^{\frac{2d+1}{4}}D_\m^{-1/2}(1+O(h)).
\ee

Let us now prove $i)$. By definition, for all $\m\in\uuu^{(0)},\ \<\phii_\m,\phii_\m\>=1$. Computing  $\<\phii_{\widehat{\m}},\phii_{\underline{\m}}\>$, using Cauchy-Schwarz inequality and noticing that $f\geq f(\widehat{\m})$ on $\supp\phii_{\widehat{\m}}$, we have
\bes
\<\phii_{\widehat{\m}},\phii_{\underline{\m}}\>=\frac{1}{\vvert{\psi_{\underline{\m}}}}\<2e^{-(f-f(\underline{\m}))/h},\phii_{\widehat{\m}}\>_{L^2(\supp\phii_{\widehat{\m}})}=\frac{1}{\vvert{\psi_{\underline{\m}}}}O\big(e^{-(f(\widehat{\m})-f(\underline{\m}))/h}\big),
\ees
which implies $i)$ using \eqref{eclpsi} and recalling $\underline{\m}$ is the lone global minimum of $f$.

Recalling $P\psi_{\underline{\m}}=P^*\psi_{\underline{\m}}=0$ we only have to prove $ii),\, iii)$ and $iv)$ for $\m=\widehat{\m}$. Using the calculus done after the poof of \cite[Proposition 5.1 $i)$]{BoLePMi22} of $\<P\psi_{\widehat{\m}},\psi_{\widehat{\m}}\>$ we have 
\bes
\<P\psi_{\widehat{\m}},\psi_{\widehat{\m}}\> = h^2C_h^{-2}\int_{\ccc_{4\tau,4\delta}}\theta^2\zeta(\ell/\tau)^2A\nabla\ell\cdot\nabla\ell \,e^{-2\big(f+\frac{\ell^2}{2}-f({\widehat{\m}})\big)/h}+O(e^{-2(\tilde S({\widehat{\m}})+\delta)/h}).
\ees

According to how we built $\ell$ and Lemma \ref{lem:xiortho},
\bes
\Big(f+\frac{\ell_0^2}{2}\Big)(\s)=f(\s),\quad \nabla\Big(f+\frac{\ell_0^2}{2}\Big)(\s)=0\ \mbox{ and }\Hess_\s\Big(f+\frac{\ell_0^2}{2}\Big)>0.
\ees
With another Laplace method, considering $-2\big(f+\frac{\ell_0^2}{2}-f({\widehat{\m}})\big)$ as the phase function and $\s$ as the global minima, we have
\bes
\<P\psi_{\widehat{\m}},\psi_{\widehat{\m}}\> = \frac{2h}{\pi}(\pi h)^{\frac{2d+1}{2}}(A\nabla\ell_0\cdot\nabla\ell_0)(\s)\Big(\det\Hess_\s\Big(f+\frac{\ell_0^2}{2}\Big)\Big)^{-1/2}e^{-2\tilde S(\widehat{\m})/h}(1+O(h)).
\ees
Using that
\bes
(A\nabla\ell_0\cdot\nabla\ell_0)(\s) = \mu = \mu(\s) \quad \mbox{and} \quad \det\Hess_\s\Big(f+\frac{\ell_0^2}{2}\Big) = D_\s^2
\ees
thanks to \eqref{eq:Axiximu} and \eqref{eq:detHplusPixi}, we have
\bes
\<P\psi_{\widehat{\m}},\psi_{\widehat{\m}}\> = \frac{2h}{\pi}(\pi h)^{\frac{2d+1}{2}}\mu D_\s^{-1} e^{-2\tilde S(\widehat{\m})/h}(1+O(h))
\ees
and \eqref{eclpsi} is enough to conclude.

Let us now prove $iii)$. Using the computations in the proof of \cite[Proposition 5.1 $iii)$]{BoLePMi22} on $\vvert{P\psi_{\widehat{\m}}}^2$, we have
\bes
\vvert{P\psi_{\widehat{\m}}}^2 = \vvert{P(\chi_\ell e^{-(f-f(\widehat{\m}))/h})}^2_{L^2(\supp\theta)} + O(e^{-2(\tilde S(\widehat{\m})+\delta)/h}),
\ees
and we recall \eqref{eq:Pchi} with the constant $C_h$
\bes
P(\chi_\ell e^{-(f-f(\widehat{\m}))/h}) = \sqrt hO(X^4+hX^2+h^2)e^{-\big(f-f(\widehat{\m})+\frac{\ell^2}{2}\big)/h}
\ees
and thus, since we are on $\supp\theta\subset\ccc_{4\tau,4\delta}$, we obtain with another Laplace method
\bes
\vvert{P\psi_{\widehat{\m}}}^2 = O(h^4)\<P\psi_{\widehat{\m}},\psi_{\widehat{\m}}\>.
\ees

Now let us move to $iv)$. We notice $P^*$ satisfies \eqref{eq:formP} with $-b$ instead of $b$ therefore, it satisfies an equation similar to \eqref{eq:Pwr}, with a $w^*$ slightly different from $w$ but with $w^*_{0,0} = w_{0,0} = 0$ and $w^*_{0,1} = O(X)$, it leads to
\bes
P^*(\chi_\ell e^{-(f-f(\m))/h})=\sqrt hO(X)e^{-\big(f-f(\m)+\frac{\ell^2}{2}\big)/h},
\ees
and hence
\bes
\vvert{P^*\psi_{\widehat{\m}}}^2 = O(h)\<P\psi_{\widehat{\m}},\psi_{\widehat{\m}}\>.
\ees

\ep

\subsection{Proof of Theorem \ref{thm3}.}
From Proposition \ref{prop:interactionmatrix} we denote
\bes
\tilde{\lambda}=\<P\phii_{\widehat{\m}},\phii_{\widehat{\m}}\>=he^{-2\tilde S(\widehat{\m})/h}\frac{\mu(\s)}{2\pi}\frac{D_{\widehat{\m}}}{D_\s}(1+O(h)).
\ees

Recalling $f(x,v,y)=\frac12\big(V(x)+\frac{|v|^2+y^2}{2}\big)$ we have that \bes
2\tilde S(\widehat{\m}) = 2(f(\s) - f(\widehat{\m})) = V(\s) - V(\widehat{\m}) = S(\widehat{\m})
\ees
and
\bes
\frac{D_{\widehat{\m}}}{D_\s} = \frac{(\det \Hess_{\widehat{\m}}f)^{1/2}}{|\det \Hess_{\s}f|^{1/2}} = \frac{2^{-(2d+1)}(\det \Hess_{\widehat{\m}}V)^{1/2}}{2^{-(2d+1)}|\det \Hess_{\s}V|^{1/2}} = \frac{(\det \Hess_{\widehat{\m}}V)^{1/2}}{|\det \Hess_{\s}V|^{1/2}},
\ees
hence
\bes
\tilde \lambda = he^{-S(\widehat{\m})/h}\frac{\mu(\s)}{2\pi}\frac{(\det \Hess_{\widehat{\m}}V)^{1/2}}{|\det \Hess_{\s}V|^{1/2}}(1+O(h)).
\ees

We thus have
\bes
\vvert{P\phii_{\widehat{\m}}}=O\Big(h^2\sqrt{\tilde{\lambda}}\Big)\quad\mbox{and}\quad\vvert{P^*\phii_{\widehat{\m}}}=O\Big(\sqrt{h\tilde{\lambda}}\Big).
\ees

Let us recall $\Pi_0=\frac1{2i\pi}\int_{\p\negthinspace D}(z-P)^{-1}dz$ where $D=D(0,c_0'g(h))$, recalling $g(h)$ is defined in \eqref{eq:g}, we also denote $u_1=\Pi_0\phii_{\widehat{\m}}$ and notice that $u_0=\Pi_0\phii_{\underline{\m}}=\phii_{\underline{\m}}$.

Using Proposition \ref{prop:interactionmatrix} $i)$ and \eqref{quasiortho}, there exists $c>0$, such that for $j,k\in\{0,1\}$
\be\label{eq:orthou}
\<u_j,u_k\> = \delta_{j,k} + O(e^{-c/h}).
\ee
Moreover, using \eqref{quasiortho}, we observe that
\bes
\begin{split}
    \<Pu_1,u_1\> &= \<P\phii_{\widehat{\m}},\phii_{\widehat{\m}}\> + \<P(\Pi_0\phii_{\widehat{\m}}-\phii_{\widehat{\m}}),\phii_{\widehat{\m}}\> + \<P\Pi_0\phii_{\widehat{\m}},\Pi_0\phii_{\widehat{\m}}-\phii_{\widehat{\m}}\>\\
    &= \tilde{\lambda} + O\big(\vvert{\Pi_0\phii_{\widehat{\m}}-\phii_{\widehat{\m}}}\vvert{P^*\phii_{\widehat{\m}}} + \vvert{\Pi_0\phii_{\widehat{\m}}-\phii_{\widehat{\m}}}\vvert{P\phii_{\widehat{\m}}}\big)\\
    &= \tilde\lambda + g(h)^{-1}\vvert{P\phii_{\widehat{\m}}}O\big(\vvert{P^*\phii_{\widehat{\m}}} + \vvert{P\phii_{\widehat{\m}}}\big)\\
    &= \tilde{\lambda} + O\big(g(h)^{-1}h^2(\sqrt h+h^2)\tilde{\lambda}\big)\\
    &= \tilde{\lambda}\big(1 + O\big(g(h)^{-1}h^\frac52\big)\big)\\
    &= \tilde{\lambda}\big(1 + O\big(\sqrt h\big)\big).
\end{split}
\ees

We then see that for $h$ small enough, $u = (u_0,u_1)$ is a basis of $\Ran\Pi_0$, thus we can consider $e_0 = u_0$, $\tilde e_1 = u_1 - \<u_1,e_0\>e_0$ and $e_1 = \frac{\tilde e_1}{\vvert{\tilde e_1}}$ (the Gram-Schmidt orthonormalization process). We have from \eqref{eq:orthou} that $\vvert{\tilde e_1} = 1 + O(e^{-c/h})$ for some $c>0$, therefore $e_1 = u_1(1 + O(e^{-c/h})) + O(e^{-c/h})e_0$. Using now that $Pe_0 = P^*e_0 = 0$, we can easily compute the matrix of $P_{|\Ran\Pi_0}$ in the basis $e = (e_0,e_1)$
\bes
\text{Mat}_eP_{|\Ran\Pi_0} = (\<Pu_{j},u_{k}\>)_{0\leq j,k\leq1}(1+O(e^{-c/h}))
\ees
and
\bes
(\<Pu_{j},u_{k}\>)_{0\leq j,k\leq1} = \begin{pmatrix}0&0\\0&\<Pu_1,u_1\>\end{pmatrix} = \begin{pmatrix}0&0\\0&\tilde{\lambda}(1+O(\sqrt h))\end{pmatrix}.
\ees
Seeing its eigenvalues, we conclude the proof of Theorem \ref{thm3}.

\newpage
\appendix
\section{Some technical results}

\subsection{Proof of Proposition \ref{prop:accretive}.}
The idea is to mimic the proof of \cite[Theorem 15.1]{He13}.

Let $h,\nu,\gamma>0$ be fixed. To show that $P$ admits a maximal accretive extension, it is first necessary to show that it is accretive, this comes from the skew-adjointness of $H_0$ and $Y$, as well as from the positivity of $\ooo$. It therefore remains to show the maximal side, for that we use the criterion which tells us that the closure of $P$ is maximal accretive if $T=P+\gamma(h/2+1)\Id$ has a dense image.

Let $f\in L^2(\R^{2d+1})$ such that
\be\label{eq:fTu}
\forall u\in\ccc_c^\infty(\R^{2d+1}),\ \< f,Tu\> =0.
\ee
We then must show that $f=0$. As $P$ is real, we can assume also is $f$. We split $Y= \underbrace{h((|v|^2-dh)\p_y-yv\cdot\p_v)}_{=Y_1}\underbrace{-\frac{y}{2}dh}_{=Y_0}$ where we can see $Y_1$ is a homogeneous differential operator of order 1 and $Y_0$ is a mere $\ccc^\infty$ function. We want to apply the standard hypoellipticity theorem for Hörmander operators to use our solution as test function, so let us verify the hypothesis : $P=-\sum_{j=1}^kX_j^2+X_0+a(x,v,y)$ with $k=d$, $X_j=\sqrt{\gamma}h\p_{v_j}$, $a(x,v,y)=\gamma\frac{|v|^2}{4}-\gamma \frac{dh}{2}-\nu \frac{y}{2}dh$ and
\bes
X_0=H_0+\nu Y_1 = v\cdot h\p_x-\p_x V\cdot h\p_v+\nu h((|v|^2-dh)\p_y-yv\cdot\p_v).
\ees
Therefore, the Lie brackets are
\bes
\begin{split}
    [X_j,X_0] &= \sqrt{\gamma}h^2\p_{x_j}+\sqrt{\gamma}\nu h^2(2v_j\p_y-y\p_{v_j}),\\
    [X_j,[X_j,X_0]] &= 2\gamma\nu h^3\p_y,
\end{split}
\ees
which ensures that $(X_j)_{j\geq0}$ is bracketgenerating. This allows us to apply Hörmander's theorem, and therefore $P$ is hypoelliptic. Because \eqref{eq:fTu} means that $Tf = 0$ in $\ddd'$, we have that $f$ is smooth, thus we can take $u=\chi f$ in \eqref{eq:fTu} for any $\chi\in\ccc^\infty_c(\R^{2d+1})$.

By taking back the computation of \cite[p219]{He13} with a modified $\zeta$ : $\zeta_k(x,v,y)=\widetilde{\zeta}\left(\frac{x}{k_1}\right)\widehat{\zeta}\left(\frac{|v|^2+y^2}{k_2}\right)$ (where $\widetilde{\zeta},\widehat{\zeta}\in\ccc_c^\infty$, are cutoffs around $0$ and $k=(k_1,k_2)\in(\R_+^*)^2$), we get
\be\label{zeta}
\begin{split}
    &\gamma h^2\vvert{\p_v(\zeta_kf)}^2+\frac{\gamma}{4}\vvert{\zeta_kvf}^2+\gamma\vvert{\zeta_kf}^2\\
    &\phantom{*****}+\int\zeta_kf^2H_0(\zeta_k)+\nu\int fY(\zeta_k^2 f) = \gamma h^2\vvert{\,|\p_v\zeta_k|f}^2.
\end{split}
\ee
Helffer does not use the expression of $\zeta_k$ before this, the result is true for any compactly supported function. For the following computations, $C$ will denote a positive constant that might change from line to line.

Noting $C(k_1) = \sup_{|x|\leq2k_1}|\nabla V|$ we have :
\bes
\begin{split}
    -\int\zeta_kf^2H_0(\zeta_k)&\leq \frac{C}{k_1}\vvert{\zeta_kvf}\vvert{f}+C\frac{C(k_1)}{k_2}\vvert{\zeta_kvf}\vvert{f}\\
    &\leq \frac{\gamma}{16}\vvert{\zeta_kvf}^2+\frac{C}{k_1^2}\vvert{f}^2+\frac{\gamma}{16}\vvert{\zeta_kvf}^2+C\frac{C(k_1)^2}{k_2^2}\vvert{f}^2.
\end{split}
\ees
And since $|\p_{v_j}\zeta_k| \leq C\frac{v_j}{k_2}$, we have $\gamma h^2\vvert{\,|\nabla_v\zeta_k|f}\leq\frac{C}{\sqrt{k_2}}\vvert{f}$. Plugging these two into \eqref{zeta} we get $$\gamma\vvert{\zeta_kf}^2+\frac{\gamma}{8}\vvert{\zeta_kvf}^2+\nu\int fY(\zeta_k^2 f) \leq C\left(\frac{1}{k_1^2}+\frac{C(k_1)^2}{k_2^2}+\frac{1}{k_2}\right)\vvert{f}^2.$$
It only remains to study the last integral, we have by integration by parts : 
\bes
\begin{array}{c}
    \int fY_1(\zeta_k^2f) = \int f(Y_1(\zeta_k)\zeta_kf+Y_1(\zeta_kf)\zeta_k)\\
    \int\zeta_kfY_1(\zeta_kf) = \frac{1}{2}\int Y_1((\zeta_kf)^2)= -\frac{h}{2}\int yv\cdot\p_v((\zeta_kf)^2)\\
    \int fY(\zeta_k^2f) = \int f^2\zeta_kY(\zeta_k)+\frac{dh}{2}\int y(\zeta_kf)^2.
\end{array}
\ees
Now with the expression of $\zeta_k$, we have $Y(\zeta_k) = -dh\left(h\p_y+\frac{y}{2}\right)(\zeta_k)$, and thus :
\bes
-\int fY(\zeta_k^2f) = h^2d\int\zeta_k f^2\p_y\zeta_k \leq \frac{C}{\sqrt{k_2}}\vvert{\zeta_kf}\vvert{f} \leq \frac{\gamma}{2\nu}\vvert{\zeta_kf}^2+\frac{C}{k_2}\vvert{f}^2.
\ees
We finally get
\bes
\frac{\gamma}{2}\vvert{\zeta_kf}^2\leq\frac{\gamma}{8}\vvert{\zeta_kvf}^2+\frac{\gamma}{2}\vvert{\zeta_kf}^2 \leq C\left(\frac{1}{k_1^2}+\frac{C(k_1)^2}{k_2^2}+\frac{1}{k_2}\right)\vvert{f}^2.
\ees
And by taking the limit $k_2\to+\infty$ and then $k_1\to+\infty$, we obtain $\vvert{f}=0$.

\ep

\subsection{Proof of Lemma \ref{lem:computprod}.}
By simple computations, we recall \eqref{eq:H_0Ydelta} :
\bes
\begin{split}
    &H_0 = v\cdot\delta_x - \p_xV\cdot\delta_v,\\
    &Y = (|v|^2-dh)\delta_y - yv\cdot\delta_v.
\end{split}
\ees
For \eqref{prod1}, this leads to
\bes
\begin{split}
    H_i^*H_j\Pr &= -H_iv_j\delta_{x_j}\Pr = (-v_iv_j\delta_{x_i}\delta_{x_j}+h\p_{x_i}V\p_{v_i}(v_j)\delta_{x_j})\Pr\\
    &= (-v_iv_j\delta_{x_i}\delta_{x_j}+\delta_{i,j}h\p_{x_i}V\delta_{x_j})\Pr.
\end{split}
\ees
For \eqref{prod2}, recalling that $Y$ is skew-adjoint and thanks to $-\delta_y = \delta_y^*-y$, we first get
\bes
\begin{split}
    Y^* &= (|v|^2-dh)(\delta_y^*-y) + yv\cdot\delta_v\\
    &= (|v|^2-dh)\delta_y^* - y(|v|^2-dh) + yv\cdot\delta_v.
\end{split}
\ees
We then get
\bes
\begin{split}
    Y^*Y\Pr &= Y^*(|v|^2-dh)\delta_y\Pr\\
    &= ((|v|^2-dh)^2\delta_y^*\delta_y - ((|v|^2-dh)^2 - v\cdot h\p_v(|v|^2-dh))y\delta_y)\Pr\\
    &= ((|v|^2-dh)^2\delta_y^*\delta_y - ((|v|^2-dh)^2 - 2h|v|^2)y\delta_y)\Pr.
\end{split}
\ees
Obtaining the last two equations \eqref{prod3} and \eqref{prod4} is pretty straightforward, we just write
\bes
\begin{split}
    H_i^*Y\Pr &= -H_i(|v|^2-dh)\delta_y\Pr\\
    &= (-v_i(|v|^2-dh)\delta_{x_i}\delta_y + \p_{x_i}Vh\p_{v_i}(|v|^2-dh)\delta_y)\Pr\\
    &= (-v_i(|v|^2-dh)\delta_{x_i}\delta_y + 2v_ih\p_{x_i}V\delta_y)\Pr,
\end{split}
\ees
and
\bes
\begin{split}
    Y^*H_i\Pr &= -Yv_i\delta_{x_i}\Pr\\
    &= (-v_i(|v|^2-dh)\delta_{x_i}\delta_y + yv\cdot h\p_v(v_i)\delta_{x_i})\Pr\\
    &= (-v_i(|v|^2-dh)\delta_{x_i}\delta_y + v_ihy\delta_{x_i})\Pr.
\end{split}
\ees

\ep

\subsection{Some resolvent estimates.}
\begin{lemma}\label{lem:bornres}
One has the following estimates
$$
\Vert(h\alpha+d\Delta_{\frac{V}{2}})(h\alpha+dB)^{-1}\Vert\leq 1
$$
and 
$$
\Vert(h\alpha+2d\nu^2hN_y)(h\alpha+dB)^{-1}\Vert\leq 1.
$$
\end{lemma}
\bp
For any $u\in L^2$, since $N_y\geq 0$, one has
\bes
\begin{split}
\Vert(h\alpha+d\Delta_{\frac{V}{2}})^{\frac 12}(h\alpha+dB)^{-\frac 12}u\Vert^2&=
\<(h\alpha+d\Delta_{\frac{V}{2}})(h\alpha+dB)^{-\frac 12}u,(h\alpha+dB)^{-\frac 12}u\>\\
&\leq \<(h\alpha+dB)(h\alpha+dB)^{-\frac 12}u,(h\alpha+dB)^{-\frac 12}u\>\\
&\leq \Vert u\Vert^2.
\end{split}
\ees
Consequently, denoting $T=(h\alpha+d\Delta_{\frac{V}{2}})^{\frac 12}(h\alpha+B)^{-\frac 12}$, one has  $\Vert T\Vert_{L^2\rightarrow L^2}\leq 1$. Moreover, since $\Delta_{\frac{V}{2}}$ and $N_y$ commute, one has
$(h\alpha+d\Delta_{\frac{V}{2}})(h\alpha+B)^{-1}=T^*T$ and  the first estimate follows immediatly. The second is identical.

\ep

\begin{lemma}\label{lem:resolvestim} There exists $C>0$ such that for all $s>0$, one has the following resolvent estimates
$$
\Vert (h\alpha+d\Delta_{\frac{V}{2}})^{-s}\Vert + \Vert (h\alpha+2d\nu^2hN_y)^{-s}\Vert\leq Ch^{-s}\alpha^{-s}
$$
and hence
$$
\Vert \Delta_{\frac{V}{2}}(h\alpha+d\Delta_{\frac{V}{2}})^{-1}\Vert + \nu^2h\Vert N_y(h\alpha+2d\nu^2hN_y)^{-1}\Vert\leq C.
$$
\end{lemma}

\bp
The first equation is straightforward when noticing $\Delta_{\frac V2}$ and $N_y$ are non negative. The second equation can be proved using functional calculus considering $x\mapsto \frac{x}{h\alpha+cx}$. Let us remark that we can compute the constant, being $2$ for the first equation, and $\frac{3}{2d}$ for the second one, but since the exact form of these constant will not be useful, we will keep $C$.

\ep

\begin{lemma}\label{lem:elliptreg} There exists $C>0$ such that for all $h>0$ and $\nu>0$ one has
$$
\Vert\delta_{x_i}(h\alpha+d\Delta_{\frac{V}{2}})^{-\frac 12}\Vert + \nu \sqrt h\Vert\delta_y(h\alpha+2d\nu^2hN_y)^{-\frac 12}\Vert\leq C.
$$
\end{lemma}
\bp
This estimate is obtained by taking the adjoint and using the spectral theorem. More precisely, for any $u\in L^2$, one has
\bes
\begin{split}
\Vert\delta_{x_i}(h\alpha+d\Delta_{\frac{V}{2}})^{-\frac 12}u\Vert^2&=\<\delta_{x_i}^*\delta_{x_i}(h\alpha+d\Delta_{\frac{V}{2}})^{-\frac 12}u, (h\alpha+d\Delta_{\frac{V}{2}})^{-\frac 12}u\>\\
&\leq \<\Delta_{\frac{V}{2}}(h\alpha+d\Delta_{\frac{V}{2}})^{-\frac 12}u, (h\alpha+d\Delta_{\frac{V}{2}})^{-\frac 12}u\>\\
&\leq \Vert \Delta_{\frac{V}{2}}(h\alpha+d\Delta_{\frac{V}{2}})^{-1}\Vert \Vert u\Vert^2\\
&\leq C \Vert u\Vert^2
\end{split}
\ees
by the previous lemma. The same arguments give the estimate on $\delta_y(h\alpha+2d\nu^2hN_y)^{-\frac 12}$.

\ep

\vfill
\bibliographystyle{siam}

\end{document}